\documentclass[10pt]{amsart}
\usepackage[T1]{fontenc}
\usepackage{geometry}
\usepackage[latin1] {inputenc}
\usepackage{amsmath}
\usepackage{amsfonts, amssymb, textcomp}
\usepackage[colorlinks=flase, linkcolor=red,urlcolor=green, citecolor=blue]{hyperref}
\usepackage{subeqnarray}

\usepackage{latexsym}
\usepackage{fancyhdr}
\usepackage{longtable}
\usepackage{amsmath, amssymb}
\usepackage{graphicx}
\usepackage{enumerate}

\numberwithin{equation}{section}

\theoremstyle{plain}
\newtheorem{proposition}{Proposition}[section]
\newtheorem{theorem}[proposition]{Theorem}
\newtheorem{lemma}[proposition]{Lemma}
\newtheorem{corollary}[proposition]{Corollary}

\newtheorem{remark}[proposition]{Remark}

\newcommand{\RR}{\mathbb{R}}

\newcommand{\NN}{\mathbb{N}}

\let\on=\operatorname

\newenvironment{demo}[1]{\par\smallskip\noindent{\bf #1.}}{\par\smallskip}

\makeatletter
\@namedef{subjclassname@2020}{%
  \textup{2020} Mathematics Subject Classification}
\makeatother

\newsavebox{\fmbox}
\newenvironment{fmpage}[1]
 {\begin{lrbox}{\fmbox}\begin{minipage}{#1}}
 {\end{minipage}\end{lrbox}\fbox{\usebox{\fmbox}}}

\title[On the equivalence between moderate growth-type conditions]
{On the equivalence between moderate growth-type conditions in the weight matrix setting}

\author[G.~Schindl]{Gerhard Schindl}

\address{G.~Schindl: Fakult\"at f\"ur Mathematik, Universit\"at Wien, Oskar-Morgenstern-Platz~1, A-1090 Wien, Austria.}
\email{gerhard.schindl@univie.ac.at}

\begin{document}

\begin{abstract}
We study the generalizations of the known equivalent reformulations of condition moderate growth from the single weight sequence to the weight matrix setting. This condition, also known in the literature under the name stability under ultradifferentiable operators, plays a significant role in the theory of ultradifferentiable (and ultraholomorophic) function classes defined in terms of weight sequences and its generalization becomes relevant when dealing with classes defined by weight matrices. In the matrix setting, we prove that the different mixed conditions are in general not equivalently satisfied anymore and we focus on weight matrices associated with (associated) weight functions.
\end{abstract}

\thanks{G. Schindl is supported by FWF-Project P33417-N}
\date{\today}

\maketitle

\section{Introduction}
In the theory of ultradifferentiable (and ultraholomorphic) function spaces defined by means of weight sequences $M=(M_p)_p$ several basic growth and regularity assumptions on the weight $M$ are arising frequently. One of the ''most prominent and classical'' conditions is {\itshape moderate growth}, denoted by \hypertarget{mg}{$(\on{mg})$}, which reads as follows:
$$\exists\;C\ge 1\;\forall\;p,q\in\NN:\;M_{p+q}\le C^{p+q} M_p M_q.$$
In \cite{Komatsu73} and in other places in the literature this requirement is denoted by $(M.2)$ and also known under the name {\itshape stability under ultradifferential operators.}\vspace{6pt}

Assumption \hyperlink{mg}{$(\on{mg})$} for $M$ implies or even characterizes important and desirable properties for the corresponding ultradifferentiable function classes of {\itshape Roumieu type} $\mathcal{E}_{\{M\}}$ and of {\itshape Beurling type} $\mathcal{E}_{(M)}$, e.g. having {\itshape Cartesian closedness} as shown in \cite{KMRc}.\vspace{6pt}

In order to apply \hyperlink{mg}{$(\on{mg})$} in proofs, several equivalent and useful variants of this condition have been shown, see \cite{PetzscheVogt}, \cite{matsumoto}, \cite{matsumotopseudo} and \cite{whitneyextensionweightmatrix}. We summarize everything in detail in Theorem \ref{mgreformulated}, for the moment we recall that \hyperlink{mg}{$(\on{mg})$} is equivalent to having
\begin{equation}\label{mgstrange}
\exists\;A\ge 1\;\forall\;p\in\NN_{>0}:\;\;\;\mu_p\le A(M_p)^{1/p},
\end{equation}
with $\mu_p:=\frac{M_p}{M_{p-1}}$. Thus \eqref{mgstrange} admits the possibility of comparing the sequences of quotients and roots of a weight sequence $M$.\vspace{6pt}

The second classical approach to introduce ultradifferentiable classes is by using a weight function $\omega:[0,+\infty)\rightarrow[0,+\infty)$, see \cite{Bjorck66} and \cite{BraunMeiseTaylor90}. Also in this setting basic growth and regularity assumptions on the weight $\omega$ are required in order to study the corresponding classes $\mathcal{E}_{\{\omega\}}$ and $\mathcal{E}_{(\omega)}$.

In \cite{BonetMeiseMelikhov07} both methods have been compared and it has been shown that in general both approaches are mutually distinct. For the case when the settings coincide \hyperlink{mg}{$(\on{mg})$} has become relevant. However, in general it is natural to ask the following: When a statement is valid for one setting, can we prove the analogous version in the other setting as well and (how) can the proofs be transferred?\vspace{6pt}

Inspired by the results shown in \cite{BonetMeiseMelikhov07}, in \cite{dissertation} and \cite{compositionpaper} ultradifferentiable classes $\mathcal{E}_{\{\mathcal{M}\}}$ and $\mathcal{E}_{(\mathcal{M})}$ defined in terms of weight matrices $\mathcal{M}=\{M^{(x)}: x\in\mathcal{I}\}$, $\mathcal{I}=\RR_{>0}$ denoting the index set, have been introduced. It has been shown that, by using the weight matrix $\mathcal{M}_{\omega}=\{W^{(x)}: x>0\}$ {\itshape associated with a given weight function} $\omega$, in this general framework one is able to treat both classical settings in a uniform and convenient way but also more classes. For the contents of this article we will not need the precise definitions of ultradifferentiable classes, we refer to \cite{compositionpaper} for more details.

By using this method one is able to see the role of moderate growth in the weight function world in a more precise way: Assuming \hyperlink{mg}{$(\on{mg})$} for some/any $W^{(x)}$ is a too strong requirement because in this case the class defined by $\mathcal{M}_{\omega}$ can already be represented by a single weight sequence class defined in terms of some arbitrary $W^{(x)}$. (Note: For the coincidence with the class defined by $\omega$, \eqref{newexpabsorb} is indispensable in order to absorb exponential growth but this property is valid by the standard assumption $\omega(2t)=O(\omega(t))$, $t\rightarrow+\infty$.)

However, by exploiting the definition of the sequences $W^{(x)}$, which is based on the weight function approach derived in \cite{BraunMeiseTaylor90}, and the convexity of the Young-conjugate $\varphi^{*}_{\omega}$ of the function $\varphi_{\omega}:=\omega\circ\exp$, one is able to see that automatically the generalization of \hyperlink{mg}{$(\on{mg})$} to a mixed condition is valid, see \eqref{newmoderategrowth}. And in general, since the matrix setting involves a parameter/index, mixed versions of the conditions arising in the weight sequence approach seem to be more natural and sufficient and thus the direct generalizations of \hyperlink{mg}{$(\on{mg})$} read as follows, see \cite[Section 4.1]{compositionpaper} and \cite[Section 7.2]{dissertation}:\par\vskip.3cm

\hypertarget{R-mg}{$(\mathcal{M}_{\{\text{mg}\}})$} \hskip1cm $\forall\;x\in\mathcal{I}\;\exists\;C>0\;\exists\;y\in\mathcal{I}\;\forall\;p,q\in\NN: M^{(x)}_{p+q}\le C^{p+q} M^{(y)}_p M^{(y)}_q$,\par\vskip.3cm

\hypertarget{B-mg}{$(\mathcal{M}_{(\text{mg})})$} \hskip1cm $\forall\;x\in\mathcal{I}\;\exists\;C>0\;\exists\;y\in\mathcal{I}\;\forall\;p,q\in\NN: M^{(y)}_{p+q}\le C^{p+q} M^{(x)}_p M^{(x)}_q$.\par\vskip.3cm

Note that in the matrix setting naturally the conditions are arising pair-wise, one has to treat conditions of Roumieu and of Beurling type, see e.g. \cite[Sect. 4.1]{compositionpaper}.\vspace{6pt}

Unfortunately, \hyperlink{R-mg}{$(\mathcal{M}_{\{\on{mg}\}})$} and \hyperlink{B-mg}{$(\mathcal{M}_{(\on{mg})})$} are not ''sufficient enough'' in order to transfer known proofs and techniques from the weight sequence setting since also knowledge on the mixed versions of the other conditions listed in Theorem \ref{mgreformulated} is relevant. Of course, a full generalization of Theorem \ref{mgreformulated} to the matrix setting is desirable, but in our attempts we have only been able to prove some partial results, see \cite[Prop. 3.6]{testfunctioncharacterization}, \cite[Thm. 9.5.1, Thm. 9.5.3]{dissertation} and \cite[Lemma 2.6]{whitneyextensionweightmatrix}. More precisely, it is not clear how the generalizations of \eqref{mgstrange}, see \eqref{rstrange} and \eqref{bstrange}, are related to \hyperlink{R-mg}{$(\mathcal{M}_{\{\on{mg}\}})$} and \hyperlink{B-mg}{$(\mathcal{M}_{(\on{mg})})$} and this loss of information led to the definition of an ''admissible weight matrix'' given in \cite[Def. 4.6]{whitneyextensionweightmatrix}.\vspace{6pt}

The aim of this article is to study and compare these mixed conditions in detail and to prove the failure of Theorem \ref{mgreformulated} in the matrix setting. Thus our results will underline the difference between the weight matrix and single weight sequence setting. More precisely, we show that the mixed version of \eqref{mgstrange}, i.e. condition \eqref{rstrange} resp. \eqref{bstrange}, is in general violated.

In order to do so we focus on weight matrices $\mathcal{M}_{\omega_M}$, where $\omega_M$ denotes the so-called {\itshape associated weight function} of a weight sequence $M$, see Section \ref{assofctsect}. The advantage in this case is that all relevant information on the matrix $\mathcal{M}_{\omega_M}$ can be expressed in terms of given $M$ in a compact form. Hence we are able to construct a counter-example in Theorem \ref{counter1} and to characterize the situation when the generalization of Theorem \ref{mgreformulated} for $\mathcal{M}_{\omega_M}$ is valid, see Proposition \ref{strangechar}. (By using the matrix $\mathcal{M}_{\omega}$ we can transfer this knowledge also to the weight function case, see Proposition \ref{strangecharom1}.)\vspace{6pt}

For the sake of completeness let us mention that in the very recent paper \cite{mixedgrowthindex} it has been shown that conditions \hyperlink{R-mg}{$(\mathcal{M}_{\{\on{mg}\}})$} and \hyperlink{B-mg}{$(\mathcal{M}_{(\on{mg})})$} can be expressed equivalently by involving new mixed growth indices related to the concept of mixed O-regular variation and by a condition involving the ''multi-index weight matrix construction'' from \cite[Sect. 9.3]{dissertation} and \cite{testfunctioncharacterization} (generalizing the situation for \hyperlink{mg}{$(\on{mg})$} which has been studied in \cite{index}).\vspace{6pt}

We summarize now the content of this paper.

After collecting and recalling all necessary definitions and assumptions on weights in Section \ref{weightscond}, in Section \ref{modgrowthsection} we are gathering and extending the knowledge on the parts of Theorem \ref{mgreformulated} which can be generalized to the matrix setting, see Theorem \ref{matrixmgtheorem}. For this we have to introduce an auxiliary matrix, denoted by $\mathcal{M}^{\on{mg}}$, but for the weight matrix $\mathcal{M}_{\omega}$ this construction is superfluous and we obtain a complete direct characterization, see Corollary \ref{matrixmgtheoremcor}.

Section \ref{failuresection} is dedicated to the study of the generalization of \eqref{mgstrange} and the failure of Theorem \ref{mgreformulated}, see Prop. \ref{strangechar} and \ref{strangecharom1} and the counter-example constructed in Section \ref{countersection}. Finally, in the Appendix \ref{snqappendix} some new formulas for the quotients of the sequences $W^{(x)}$, which are needed in the proof of the main statement, are applied to the $\liminf$-conditions arising in \cite[Thm. 3.11 $(v)$]{index} (and generalizing the {\itshape strong non-quasianalyticity condition} for weight sequences), see Proposition \ref{strongnonquasilemma}.

\section{Weights and conditions}\label{weightscond}

\subsection{General notation}
We write $\NN:=\{0,1,2,\dots\}$ and $\NN_{>0}:=\{1,2,3,\dots\}$.

\subsection{Weight sequences}\label{weightsequences}
Given a sequence $M=(M_p)_p\in\RR_{>0}^{\NN}$ we also use the notation $\mu=(\mu_p)_p$ with $\mu_p:=\frac{M_p}{M_{p-1}}$, $p\ge 1$, $\mu_0:=1$, and analogously for all other arising sequences. $M$ is called {\itshape normalized} if $1=M_0\le M_1$ holds true.\vspace{6pt}

$M$ is called {\itshape log-convex} if
$$\forall\;p\in\NN_{>0}:\;M_p^2\le M_{p-1} M_{p+1},$$
equivalently if $\mu$ is nondecreasing. If $M$ is log-convex and normalized, then both $p\mapsto M_p$ and $p\mapsto(M_p)^{1/p}$ are nondecreasing, $(M_p)^{1/p}\le\mu_p$ for all $p\in\NN_{>0}$ and finally $M_pM_q\le M_{p+q}$ for all $p,q\in\NN$.\vspace{6pt}

For our purpose it is convenient to consider the following set of sequences
$$\hypertarget{LCset}{\mathcal{LC}}:=\{M\in\RR_{>0}^{\NN}:\;M\;\text{is normalized, log-convex},\;\lim_{p\rightarrow+\infty}(M_p)^{1/p}=+\infty\}.$$
We see that $M\in\hyperlink{LCset}{\mathcal{LC}}$ if and only if $1=\mu_0\le\mu_1\le\dots$ and $\lim_{p\rightarrow+\infty}\mu_p=+\infty$ (e.g. see \cite[p. 104]{compositionpaper}) and there is a one-to-one correspondence between $M$ and $\mu=(\mu_p)_p$ by taking $M_p:=\prod_{i=0}^p\mu_i$.\vspace{6pt}

$M$ has {\itshape derivation closedness}, denoted by \hypertarget{dc}{$(\text{dc})$}, if
$$\exists\;D\ge 1\;\forall\;p\in\NN:\;\;\;M_{p+1}\le D^{p+1} M_p\Longleftrightarrow\mu_{p+1}\le D^{p+1}.$$
In \cite{Komatsu73} this condition is denoted by $(M.2')$ and it is immediate that \hyperlink{mg}{$(\on{mg})$} implies derivation closedness.\vspace{6pt}

We say that $M$ has condition \hypertarget{beta1}{$(\beta_1)$} (introduced in \cite{petzsche}), if
$$\exists\;Q\in\NN_{\ge 2}:\;\;\;\liminf_{p\rightarrow+\infty}\frac{\mu_{Qp}}{\mu_p}>Q.$$

Moreover, there for $M\in\hyperlink{LCset}{\mathcal{LC}}$ it has been shown that \hyperlink{beta1}{$(\beta_1)$} is equivalent to requiring \hyperlink{gamma1}{$(\gamma_1)$}, i.e.
$$\sup_{p\in\NN_{>0}}\frac{\mu_p}{p}\sum_{k\ge p}\frac{1}{\mu_k}<+\infty.$$
In the literature \hyperlink{gamma1}{$(\gamma_1)$} is also called ''strong nonquasianalyticity condition'' and in \cite{Komatsu73} it is denoted by $(M.3)$ (in fact, there $\frac{\mu_p}{p}$ is replaced by $\frac{\mu_p}{p-1}$ for $p\ge 2$ but which is equivalent to having \hyperlink{gamma1}{$(\gamma_1)$}).\vspace{6pt}


A weaker requirement than \hyperlink{beta1}{$(\beta_1)$} is
\begin{equation}\label{beta3}
\exists\;Q\in\NN_{\ge 2}:\;\;\;\liminf_{p\rightarrow+\infty}\frac{\mu_{Qp}}{\mu_p}>1,
\end{equation}
which is arising in the main characterizing results in \cite{BonetMeiseMelikhov07} and denoted by $(\beta_3)$ in \cite{dissertation}. Conditions of this type are also showing up in \cite[Thm. 3.11 $(v)$]{index}.\vspace{6pt}

Let $M,N\in\RR_{>0}^{\NN}$ be given, we write $M\hypertarget{preceq}{\preceq}N$ if $\sup_{p\in\NN_{>0}}\left(\frac{M_p}{N_p}\right)^{1/p}<+\infty$ and call $M$ and $N$ {\itshape equivalent}, denoted by $M\hypertarget{approx}{\approx}N$, if $M\hyperlink{preceq}{\preceq}N$ and $N\hyperlink{preceq}{\preceq}M$. Property \hyperlink{mg}{$(\on{mg})$} is clearly preserved under \hyperlink{approx}{$\approx$} and for \hyperlink{beta1}{$(\beta_1)$} this follows by the characterizations obtained in \cite{petzsche}.

If $M\hyperlink{approx}{\approx}N$, then $\mathcal{E}_{\{M\}}=\mathcal{E}_{\{N\}}$ and $\mathcal{E}_{(M)}=\mathcal{E}_{(N)}$ as locally convex vector spaces, see e.g. \cite[Prop. 2.11 $(1)$]{compositionpaper}.

Finally, write $M\le N$ if $M_p\le N_p$ for all $p\in\NN$.\vspace{6pt}


\subsection{Associated weight function}\label{assofctsect}
Let $M\in\RR_{>0}^{\NN}$ (with $M_0=1$), then the {\itshape associated function} $\omega_M: \RR_{\ge 0}\rightarrow\RR\cup\{+\infty\}$ is defined by
\begin{equation*}\label{assofunc}
\omega_M(t):=\sup_{p\in\NN}\log\left(\frac{t^p}{M_p}\right)\;\;\;\text{for}\;t>0,\hspace{30pt}\omega_M(0):=0.
\end{equation*}
For an abstract introduction of the associated function we refer to \cite[Chapitre I]{mandelbrojtbook}, see also \cite[Definition 3.1]{Komatsu73}. If $\liminf_{p\rightarrow+\infty}(M_p)^{1/p}>0$, then $\omega_M(t)=0$ for sufficiently small $t$, since $\log\left(\frac{t^p}{M_p}\right)<0\Leftrightarrow t<(M_p)^{1/p}$ holds for all $p\in\NN_{>0}$. Moreover under this assumption $t\mapsto\omega_M(t)$ is a continuous nondecreasing function, which is convex in the variable $\log(t)$ and tends faster to infinity than any $\log(t^p)$, $p\ge 1$, as $t\rightarrow+\infty$. $\lim_{p\rightarrow+\infty}(M_p)^{1/p}=+\infty$ implies that $\omega_M(t)<+\infty$ for each finite $t$ which shall be considered as a basic assumption for defining $\omega_M$.\vspace{6pt}

By definition, the following is immediate: If $M,N\in\hyperlink{LCset}{\mathcal{LC}}$ are given such that $M\hyperlink{preceq}{\preceq}N$, then
\begin{equation}\label{assinclusion}
\exists\;C\ge 1\;\forall\;t\ge 0:\;\;\;\omega_N(t)\le\omega_M(Ct).
\end{equation}
We define the {\itshape counting function} $\Sigma_M$ by
\begin{equation}\label{counting}
\Sigma_M(t):=|\{p\in\NN_{>0}:\;\;\;\mu_p\le t\}|.
\end{equation}
It is known that for given $M\in\hyperlink{LCset}{\mathcal{LC}}$ the functions $\omega_M$ and $\Sigma_M$ are related by the following integral representation formula, see \cite[1.8. III]{mandelbrojtbook} and also \cite[$(3.11)$]{Komatsu73}:
\begin{equation}\label{intrepr}
\omega_M(t)=\int_0^t\frac{\Sigma_M(u)}{u}du=\int_{\mu_1}^t\frac{\Sigma_M(u)}{u}du.
\end{equation}
Consequently, $\omega_M$ vanishes on $[0,\mu_1]$, in particular on the unit interval.

Finally, if $M\in\hyperlink{LCset}{\mathcal{LC}}$, then we can compute $M$ by involving $\omega_M$ as follows, see \cite[Chapitre I, 1.4, 1.8]{mandelbrojtbook} and also \cite[Prop. 3.2]{Komatsu73}:
\begin{equation}\label{Prop32Komatsu}
M_p=\sup_{t\ge 0}\frac{t^p}{\exp(\omega_{M}(t))},\;\;\;p\in\NN.
\end{equation}

\subsection{Weight functions}\label{weightfunctions}
A function $\omega:[0,+\infty)\rightarrow[0,+\infty)$ is called a {\itshape weight function} (in the terminology of \cite[Section 2.1]{index} and \cite[Section 2.2]{sectorialextensions}), if it is continuous, nondecreasing, $\omega(0)=0$ and $\lim_{t\rightarrow+\infty}\omega(t)=+\infty$. If $\omega$ satisfies in addition $\omega(t)=0$ for all $t\in[0,1]$, then we call $\omega$ a {\itshape normalized weight function}. For convenience we will write that $\omega$ has $\hypertarget{om0}{(\omega_0)}$ if it is a normalized weight.\vspace{6pt}


Moreover we consider the following conditions; this list of properties has already been used in ~\cite{dissertation}.

\begin{itemize}
\item[\hypertarget{om1}{$(\omega_1)}$] $\omega(2t)=O(\omega(t))$ as $t\rightarrow+\infty$, i.e. $\exists\;L\ge 1\;\forall\;t\ge 0:\;\;\;\omega(2t)\le L(\omega(t)+1)$.


\item[\hypertarget{om3}{$(\omega_3)$}] $\log(t)=o(\omega(t))$ as $t\rightarrow+\infty$.

\item[\hypertarget{om4}{$(\omega_4)$}] $\varphi_{\omega}:t\mapsto\omega(e^t)$ is a convex function on $\RR$.


\item[\hypertarget{om6}{$(\omega_6)$}] $\exists\;H\ge 1\;\forall\;t\ge 0:\;2\omega(t)\le\omega(H t)+H$.
\end{itemize}

Finally, we recall the {\itshape strong non-quasianalyticity} condition for weight functions
\begin{equation}\label{assostrongnq}
\exists\;C\ge 1\;\forall\;y\ge 0:\;\;\;\int_1^{+\infty}\frac{\omega(yt)}{t^2}dt\le C\omega(y)+C.
\end{equation}

For convenience we define the set
$$\hypertarget{omset0}{\mathcal{W}_0}:=\{\omega:[0,\infty)\rightarrow[0,\infty): \omega\;\text{has}\;\hyperlink{om0}{(\omega_0)},\hyperlink{om3}{(\omega_3)},\hyperlink{om4}{(\omega_4)}\}.$$
For any $\omega\in\hyperlink{omset0}{\mathcal{W}_0}$ we define the {\itshape Legendre-Fenchel-Young-conjugate} of $\varphi_{\omega}$ by
\begin{equation}\label{legendreconjugate}
\varphi^{*}_{\omega}(x):=\sup\{x y-\varphi_{\omega}(y): y\ge 0\},\;\;\;x\ge 0,
\end{equation}
with the following properties, e.g. see \cite[Remark 1.3, Lemma 1.5]{BraunMeiseTaylor90}: It is convex and nondecreasing, $\varphi^{*}_{\omega}(0)=0$, $\varphi^{**}_{\omega}=\varphi_{\omega}$, $\lim_{x\rightarrow+\infty}\frac{x}{\varphi^{*}_{\omega}(x)}=0$ and finally $x\mapsto\frac{\varphi_{\omega}(x)}{x}$ and $x\mapsto\frac{\varphi^{*}_{\omega}(x)}{x}$ are nondecreasing on $[0,+\infty)$. Note that by normalization we can extend the supremum in \eqref{legendreconjugate} from $y\ge 0$ to $y\in\RR$ without changing the value of $\varphi^{*}_{\omega}(x)$ for given $x\ge 0$.

\vspace{6pt}

Let $\sigma,\tau$ be weight functions, we write $\sigma\hypertarget{ompreceq}{\preceq}\tau$ if $\tau(t)=O(\sigma(t))\;\text{as}\;t\rightarrow+\infty$
and call them equivalent, denoted by $\sigma\hypertarget{sim}{\sim}\tau$, if
$\sigma\hyperlink{ompreceq}{\preceq}\tau$ and $\tau\hyperlink{ompreceq}{\preceq}\sigma$. If $\sigma\hyperlink{sim}{\sim}\tau$, then $\mathcal{E}_{\{\sigma\}}=\mathcal{E}_{\{\tau\}}$ and $\mathcal{E}_{(\sigma)}=\mathcal{E}_{(\tau)}$ as locally convex vector spaces, see e.g. \cite[Cor. 5.17 $(1)$]{compositionpaper}.\vspace{6pt}

We recall the following known result, e.g. see \cite[Lemma 2.8]{testfunctioncharacterization} resp. \cite[Lemma 2.4]{sectorialextensions} and the references mentioned in the proofs there.

\begin{lemma}\label{assoweightomega0}
Let $M\in\hyperlink{LCset}{\mathcal{LC}}$, then $\omega_M\in\hyperlink{omset0}{\mathcal{W}_0}$ holds true and \hyperlink{om6}{$(\omega_6)$} for $\omega_M$ if and only if $M$ has \hyperlink{mg}{$(\on{mg})$}.
\end{lemma}

\subsection{Weight matrices}\label{classesweightmatrices}
For the following definitions and conditions see also \cite[Section 4]{compositionpaper}.

Let $\mathcal{I}=\RR_{>0}$ denote the index set (equipped with the natural order), a {\itshape weight matrix} $\mathcal{M}$ associated with $\mathcal{I}$ is a (one parameter) family of weight sequences $\mathcal{M}:=\{M^{(x)}\in\RR_{>0}^{\NN}: x\in\mathcal{I}\}$, such that
$$\forall\;x\in\mathcal{I}:\;M^{(x)}\;\text{is normalized, nondecreasing},\;M^{(x)}\le M^{(y)}\;\text{for}\;x\le y.$$
We call a weight matrix $\mathcal{M}$ {\itshape standard log-convex,} denoted by \hypertarget{Msc}{$(\mathcal{M}_{\on{sc}})$}, if
$$\forall\;x\in\mathcal{I}:\;M^{(x)}\in\hyperlink{LCset}{\mathcal{LC}}.$$
Moreover, we have the corresponding sequences of quotients $\mu^{(x)}$ given by $\mu^{(x)}_p:=\frac{M^{(x)}_p}{M^{(x)}_{p-1}}$ for $p\in\NN_{>0}$, $\mu^{(x)}_0:=1$.

A matrix is called {\itshape constant} if $M^{(x)}\hyperlink{approx}{\approx}M^{(y)}$ for all $x,y\in\mathcal{I}$. In particular, $\mathcal{M}=\{M\}$ is constant (when formally setting $M^{(x)}=M$ for any $x\in\mathcal{I}$).\vspace{6pt}

Let $\mathcal{M}=\{M^{(x)}: x>0\}$ and $\mathcal{N}=\{N^{(y)}: y>0\}$ be given. We write $\mathcal{M}\hypertarget{Mroumpreceq}{\{\preceq\}}\mathcal{N}$ if
$$\forall\;x>0\;\exists\;y>0:\;\;\;M^{(x)}\hyperlink{preceq}{\preceq}N^{(y)},$$
and $\mathcal{M}\hypertarget{Mbeurpreceq}{(\preceq)}\mathcal{N}$ if
$$\forall\;x>0\;\exists\;y>0:\;\;\;M^{(y)}\hyperlink{preceq}{\preceq}N^{(x)}.$$
We call $\mathcal{M}$ and $\mathcal{N}$ to be $R$-equivalent, if
$\mathcal{M}\hyperlink{Mroumpreceq}{\{\preceq\}}\mathcal{N}$ and $\mathcal{N}\hyperlink{Mroumpreceq}{\{\preceq\}}\mathcal{M}$ and $B$-equivalent, if $\mathcal{M}\hyperlink{Mbeurpreceq}{(\preceq)}\mathcal{N}$ and $\mathcal{N}\hyperlink{Mbeurpreceq}{(\preceq)}\mathcal{M}$.\vspace{6pt}

If $\mathcal{M}$ and $\mathcal{N}$ are $R$-equivalent, then $\mathcal{E}_{\{\mathcal{M}\}}=\mathcal{E}_{\{\mathcal{N}\}}$ and if they are $B$-equivalent, then $\mathcal{E}_{(\mathcal{M})}=\mathcal{E}_{(\mathcal{N})}$ as locally convex vector spaces, see e.g. \cite[Prop. 4.6 $(1)$]{compositionpaper}.

\begin{remark}\label{constantremark}
Let a non-constant matrix $\mathcal{M}=\{M^{(x)}: x\in\mathcal{I}\}$ be given. Assume that there exists $x_0\in\mathcal{I}$ such that $M^{(y)}\hyperlink{approx}{\approx}M^{(x)}$ for all $x,y\ge x_0$. Then, when dealing with Roumieu type classes, we can replace $\mathcal{M}$ by $\mathcal{N}=\{M^{(x_0)}\}$ since $\mathcal{E}_{\{\mathcal{M}\}}=\mathcal{E}_{\{M^{x_0}\}}$.

Similarly, when there exists $x_0\in\mathcal{I}$ such that $M^{(y)}\hyperlink{approx}{\approx}M^{(x)}$ for all $x,y\le x_0$, then in the Beurling setting we can replace $\mathcal{M}$ by $\mathcal{N}=\{M^{(x_0)}\}$ since $\mathcal{E}_{(\mathcal{M})}=\mathcal{E}_{(M^{x_0})}$.

In particular, these comments apply to the case when $\mathcal{M}$ consists of only finitely many (w.l.o.g. pair-wise non-equivalent) sequences.
\end{remark}

We summarize some facts which are shown in \cite[Sect. 5]{compositionpaper} and \cite[Sect. 4, Sect. 5]{dissertation} and are needed in this work. All properties listed below are valid for $\omega\in\hyperlink{omset0}{\mathcal{W}_0}$, except \eqref{newexpabsorb} for which \hyperlink{om1}{$(\omega_1)$} is necessary (and underlines the importance of this condition in this context).

\begin{itemize}
\item[$(i)$] The idea was that to each $\omega\in\hyperlink{omset0}{\mathcal{W}_0}$ we can associate a standard log-convex weight matrix $\mathcal{M}_{\omega}:=\{W^{(l)}=(W^{(l)}_p)_{p\in\NN}: l>0\}$ by\vspace{6pt}

    \centerline{$W^{(l)}_p:=\exp\left(\frac{1}{l}\varphi^{*}_{\omega}(lp)\right)$.}\vspace{6pt}

For the corresponding sequence of quotients we write $\vartheta^{(l)}$.

\item[$(ii)$] $\mathcal{M}_{\omega}$ satisfies
    \begin{equation}\label{newmoderategrowth}
    \forall\;l>0\;\forall\;p,q\in\NN:\;\;\;W^{(l)}_{p+q}\le W^{(2l)}_pW^{(2l)}_q.
    \end{equation}
    so both \hyperlink{R-mg}{$(\mathcal{M}_{\{\on{mg}\}})$} and \hyperlink{B-mg}{$(\mathcal{M}_{(\on{mg})})$} are satisfied.

\item[$(iii)$] \hyperlink{om6}{$(\omega_6)$} holds if and only if some/each $W^{(l)}$ satisfies \hyperlink{mg}{$(\on{mg})$} if and only if $W^{(l)}\hyperlink{approx}{\approx}W^{(n)}$ for each $l,n>0$. Consequently \hyperlink{om6}{$(\omega_6)$} is characterizing the situation when $\mathcal{M}_{\omega}$ is constant.

\item[$(iv)$] In case $\omega$ has in addition \hyperlink{om1}{$(\omega_1)$}, then $\mathcal{M}_{\omega}$ has also
     \begin{equation}\label{newexpabsorb}
     \forall\;h\ge 1\;\exists\;A\ge 1\;\forall\;l>0\;\exists\;D\ge 1\;\forall\;p\in\NN:\;\;\;h^pW^{(l)}_p\le D W^{(Al)}_p,
     \end{equation}
and this estimate is crucial for proving $\mathcal{E}_{\{\mathcal{M}_{\omega}\}}=\mathcal{E}_{\{\omega\}}$ and $\mathcal{E}_{(\mathcal{M}_{\omega})}=\mathcal{E}_{(\omega)}$ as locally convex vector spaces.

\item[$(v)$] We have $\omega\hyperlink{sim}{\sim}\omega_{W^{(l)}}$ for each $l>0$, more precisely
\begin{equation}\label{goodequivalenceclassic}
\forall\;l>0\,\,\exists\,D_{l}>0\;\forall\;t\ge 0:\;\;\;l\omega_{W^{(l)}}(t)\le\omega(t)\le 2l\omega_{W^{(l)}}(t)+D_{l},
\end{equation}
see \cite[Theorem 4.0.3, Lemma 5.1.3]{dissertation} and also \cite[Lemma 2.5]{sectorialextensions}.
\end{itemize}

Since for any given $M\in\hyperlink{LCset}{\mathcal{LC}}$ we have $\omega_M\in\hyperlink{omset0}{\mathcal{W}_0}$, see Lemma \ref{assoweightomega0}, it makes sense to define the matrix associated with the weight $\omega_M$ by
$$\mathcal{M}_{\omega_M}:=\{M^{(l)}: l>0\}.$$
Then we get
\begin{equation}\label{Mxonestable}
\forall\;p\in\NN:\;\;\;M_p=M^{(1)}_p,
\end{equation}
which follows by applying \eqref{Prop32Komatsu} (see also the proof of \cite[Thm. 6.4]{testfunctioncharacterization}):
 \begin{align*}
 M^{(1)}_p&:=\exp(\varphi^{*}_{\omega_{M}}(p))=\exp(\sup_{y\ge 0}\{py-\omega_{M}(e^y)\})=\sup_{y\ge 0}\exp(py-\omega_{M}(e^y))
 \\&
 =\sup_{y\ge 0}\frac{\exp(py)}{\exp(\omega_{M}(e^y))}=\sup_{t\ge 1}\frac{t^p}{\exp(\omega_{M}(t))}=\sup_{t\ge 0}\frac{t^p}{\exp(\omega_{M}(t))}=M_p.
 \end{align*}
For this recall that we have $\omega_{M}(t)=0$ for $0\le t\le 1$, see \eqref{intrepr}.

Moreover, one has by definition
\begin{equation}\label{Mgoodtransform}
\forall\;l\in\NN_{>0}\;\forall\;p\in\NN:\;\;\;M^{(l)}_p=\exp\left(\frac{1}{l}\varphi^{*}_{\omega_{M}}(l p)\right)=(M^{(1)}_{l p})^{1/l}=(M_{l p})^{1/l}.
\end{equation}

In particular \eqref{Mgoodtransform} holds for the sequences $W^{(l)}\in\mathcal{M}_{\omega}$ (except the very last equality).

\section{Moderate growth conditions for abstractly given weight matrices}\label{modgrowthsection}
\subsection{Known characterization for the single weight setting}
For given $M\in\hyperlink{LCset}{\mathcal{LC}}$ in the literature there exist several known equivalent reformulations of \hyperlink{mg}{$(\on{mg})$}. We refer to \cite[Lemma 5.3]{PetzscheVogt}, \cite{matsumoto}, \cite[Appendix B]{matsumotopseudo} and finally to \cite[Lemma 2.2]{whitneyextensionweightmatrix} and summarize everything in the next statement (in particular extending Lemma \ref{assoweightomega0}).

\begin{theorem}\label{mgreformulated}
Let $M\in\hyperlink{LCset}{\mathcal{LC}}$ be given, then the following are equivalent:
\begin{itemize}
\item[$(i)$] $M$ has \hyperlink{mg}{$(\on{mg})$},

\item[$(ii)$] $M$ satisfies
$$\exists\;A\ge 1\;\forall\;p\in\NN:\;\;\;M_{2p}\le A^{2p}(M_p)^2,$$

\item[$(iii)$] $\mu$ satisfies $\sup_{p\in\NN}\frac{\mu_{2p}}{\mu_p}<+\infty$,

\item[$(iv)$] $\omega_M$ satisfies \hyperlink{om6}{$(\omega_6)$},

\item[$(v)$] $\Sigma_M$ satisfies \hyperlink{om6}{$(\omega_6)$},

\item[$(vi)$] $M$ has \eqref{mgstrange}.
\end{itemize}
\end{theorem}

Concerning this result we remark that:

\begin{itemize}
\item[$(a)$] \eqref{mgstrange} together with log-convexity and normalization ensure that the sequences of quotients and roots are comparable up to a constant.

\item[$(b)$] Condition \hyperlink{mg}{$(\on{mg})$} enables a special flexibility: We consider the transformation
\begin{equation}\label{pitrafo}
\pi^s:\;\;\;(M_p)_p\mapsto(p!^sM_p)_p,\;\;\;s\in\RR.
\end{equation}
By Stirling's formula it is immediate to see that assertions $(i)$, $(ii)$, $(iii)$ and $(vi)$ are preserved under $\pi^s$ for any $s\in\RR$. But note that in general for $s<0$ one will lose the assumption log-convexity and also $(p!^sM_p)^{1/p}\rightarrow+\infty$ as $p\rightarrow+\infty$ might be violated.
\end{itemize}

\subsection{The abstract weight matrix setting}
We want to see if resp. which parts of Theorem \ref{mgreformulated} can be generalized to the weight matrix setting.

First, let us summarize the following consequences for conditions \hyperlink{R-mg}{$(\mathcal{M}_{\{\on{mg}\}})$} and \hyperlink{B-mg}{$(\mathcal{M}_{(\on{mg})})$}:

\begin{itemize}
\item[$(*)$] Since $M^{(x)}\le M^{(y)}$ for $x\le y$, we can assume $y\ge x$ in \hyperlink{R-mg}{$(\mathcal{M}_{\{\on{mg}\}})$} and $y\le x$ in \hyperlink{B-mg}{$(\mathcal{M}_{(\on{mg})})$}. In both conditions we can take $x=y$ if and only if $M^{(x)}$ has \hyperlink{mg}{$(\on{mg})$}.

\item[$(*)$] If $\mathcal{M}=\{M\}$, resp. more generally if $\mathcal{M}=\{M^{(x)}: x\in\mathcal{I}\}$ is constant, then $\mathcal{M}$ has \hyperlink{R-mg}{$(\mathcal{M}_{\{\on{mg}\}})$} and/or \hyperlink{B-mg}{$(\mathcal{M}_{(\on{mg})})$} if and only if $M$ has \hyperlink{mg}{$(\on{mg})$} resp. some/each $M^{(x)}$ has \hyperlink{mg}{$(\on{mg})$}.

\item[$(*)$] If $\mathcal{M}\equiv\mathcal{M}_{\omega}$, with $\mathcal{M}_{\omega}$ being the matrix associated with a given weight $\omega\in\hyperlink{omset0}{\mathcal{W}_0}$, then by \eqref{newmoderategrowth} both \hyperlink{R-mg}{$(\mathcal{M}_{\{\on{mg}\}})$} and \hyperlink{B-mg}{$(\mathcal{M}_{(\on{mg})})$} hold true.

\item[$(*)$] Obviously, \hyperlink{R-mg}{$(\mathcal{M}_{\{\on{mg}\}})$} is preserved under $R$-equivalence and \hyperlink{B-mg}{$(\mathcal{M}_{(\on{mg})})$} under $B$-equivalence of weight matrices.
\end{itemize}

By combining \cite[Prop. 3.6]{testfunctioncharacterization} and \cite[Thm. 9.5.1, Thm. 9.5.3]{dissertation}, we recall the following known characterization in the matrix setting which gives a partial generalization of Theorem \ref{mgreformulated}.

\begin{proposition}\label{firstmg}
Let $\mathcal{M}=\{M^{(x)}: x\in\mathcal{I}\}$ be \hyperlink{Msc}{$(\mathcal{M}_{\on{sc}})$}, then in the Roumieu setting the following are equivalent:
\begin{itemize}

\item[$(I_R)$] $\mathcal{M}$ has $\hyperlink{R-mg}{(\mathcal{M}_{\{\on{mg}\}})}$,

\item[$(II_R)$] we get
$$\forall\;x\in\mathcal{I}\;\exists\;C>0\;\exists\;y\in\mathcal{I}\;\forall\;p\in\NN:\;\;\; M^{(x)}_{2p}\le C^{2p}(M^{(y)}_p)^2,$$

\item[$(III_R)$] we get
$$\forall\;x\in\mathcal{I}\;\exists\;H\ge 1\;\exists\;y\in\mathcal{I}\;\forall\;t\ge 0:\;\;\; 2\omega_{M^{(y)}}(t)\le\omega_{M^{(x)}}(Ht)+H.$$
\end{itemize}

Moreover, in the Beurling setting, we have the following equivalences:
\begin{itemize}

\item[$(I_B)$] $\mathcal{M}$ has $\hyperlink{B-mg}{(\mathcal{M}_{(\on{mg})})}$,

\item[$(II_B)$] we get
$$\forall\;y\in\mathcal{I}\;\exists\;C>0\;\exists\;x\in\mathcal{I}\;\forall\;p\in\NN:\;\;\; M^{(x)}_{2p}\le C^{2p}(M^{(y)}_p)^2,$$

\item[$(III_B)$] we get
$$\forall\;x\in\mathcal{I}\;\exists\;H\ge 1\;\exists\;y\in\mathcal{I}\;\forall\;t\ge 0:\;\;\; 2\omega_{M^{(x)}}(t)\le\omega_{M^{(y)}}(Ht)+H.$$
\end{itemize}
\end{proposition}

{\itshape Note:}

\begin{itemize}
\item[$(*)$] In particular, Proposition \ref{firstmg} applies to $\mathcal{M}_{\omega}$, $\omega\in\hyperlink{omset0}{\mathcal{W}_0}$.

\item[$(*)$] Assertion $(III_R)$ resp. $(III_B)$ in the previous result is the mixed \hyperlink{om6}{$(\omega_6)$}-condition of the particular type. But even if all associated weight functions are equivalent w.r.t. $\hyperlink{sim}{\sim}$, then in general we cannot conclude that \hyperlink{om6}{$(\omega_6)$} for each/some $\omega_{M^{(x)}}$ is following.

\item[$(*)$] It is straightforward to see (separately) that all arising assertions in the previous result are stable under $R$- resp. $B$-equivalence of weight matrices, for $(III_R)$ resp. $(III_B)$ we involve \eqref{assinclusion}.

\item[$(*)$] By Stirling's formula it is also immediate that $(I_R)$, $(II_R)$, $(I_B)$ and $(II_B)$ are preserved under the mapping $\pi^s$, see \eqref{pitrafo}.
\end{itemize}
\vspace{6pt}

Concerning the assertion $(III_R)$ resp. $(III_B)$ we can see the following result.

\begin{proposition}\label{secondmg}
Let $\mathcal{M}=\{M^{(x)}: x\in\mathcal{I}\}$ be \hyperlink{Msc}{$(\mathcal{M}_{\on{sc}})$}.

We consider in the Roumieu setting the following assertions:
\begin{itemize}

\item[$(IV_R)$] $\mathcal{M}$ satisfies
$$\forall\;x\in\mathcal{I}\;\exists\;A\ge 1\;\exists\;y\in\mathcal{I}\;\forall\;p\in\NN:\;\;\;\mu^{(x)}_{2p}\le A\mu^{(y)}_p,$$

\item[$(V_R)$] we have
$$\forall\;x\in\mathcal{I}\;\exists\;A\ge 1\;\exists\;y\in\mathcal{I}\;\forall\;t\ge 0:\;\;\;2\Sigma_{M^{(y)}}(t)\le\Sigma_{M^{(x)}}(At).$$

\end{itemize}

Moreover, in the Beurling setting, we consider the following assertions:
\begin{itemize}
\item[$(IV_B)$] $\mathcal{M}$ satisfies
$$\forall\;x\in\mathcal{I}\;\exists\;A\ge 1\;\exists\;y\in\mathcal{I}\;\forall\;p\in\NN:\;\;\;\mu^{(y)}_{2p}\le A\mu^{(x)}_p,$$

\item[$(V_B)$] we have
$$\forall\;x\in\mathcal{I}\;\exists\;A\ge 1\;\exists\;y\in\mathcal{I}\;\forall\;t\ge 0:\;\;\;2\Sigma_{M^{(x)}}(t)\le\Sigma_{M^{(y)}}(At).$$

\end{itemize}
Then we get $(IV_R)\Leftrightarrow(V_R)\Rightarrow(III_R)$ and $(IV_B)\Leftrightarrow(V_B)\Rightarrow(III_B)$.
\end{proposition}

\demo{Proof}
$(IV_R)\Leftrightarrow(V_R)$ resp. $(IV_B)\Leftrightarrow(V_B)$ follows by the definition of the counting function $\Sigma_M$ in \eqref{counting} and since each sequence is log-convex which is equivalent to the fact that each sequence of quotients $\mu^{(x)}$ is nondecreasing.\vspace{6pt}

$(V_R)\Rightarrow(III_R)$ resp. $(V_B)\Rightarrow(III_B)$ follows by involving the integral representation formula \eqref{intrepr}, see also \cite[Lemma 2.2]{whitneyextensionweightmatrix} for the single weight sequence case.

Alternatively, we can prove $(IV_R)\Rightarrow(II_R)$ resp. $(IV_B)\Rightarrow(II_B)$ analogously as in \cite[Lemma 5.3 $(1)\Rightarrow(4)$]{PetzscheVogt}. We only consider the Roumieu case in detail and estimate by log-convexity for $M^{(x)}$ as follows:
$$M^{(x)}_{2p}=\mu^{(x)}_1\cdots\mu^{(x)}_{2p}=\prod_{i=1}^{p}\mu^{(x)}_{2i-1}\mu^{(x)}_{2i}\le\prod_{i=1}^{p}(\mu^{(x)}_{2i})^2\le\prod_{i=1}^pA^2(\mu^{(y)}_i)^2=A^{2p}(M^{(y)}_p)^2.$$
\qed\enddemo

Next we want to see if all assertions $(I_R)-(V_R)$ resp. $(I_B)-(V_B)$ are equivalent or at least satisfied simultaneously. For this let now $\mathcal{M}=\{M^{(x)}: x\in\mathcal{I}\}$ be given and assume that $\mathcal{M}$ is \hyperlink{Msc}{$(\mathcal{M}_{\on{sc}})$}. Then, for technical reasons, we introduce the ''shifted matrix'' $\widetilde{\mathcal{M}}:=\{\widetilde{M}^{(x)}: x\in\mathcal{I}\}$ by setting
\begin{equation}\label{tildedef}
\widetilde{M}^{(x)}_p:=(M^{(x)}_{4p})^{1/4},\;\;\;p\in\NN.
\end{equation}
Note that $\widetilde{\mathcal{M}}$ is again \hyperlink{Msc}{$(\mathcal{M}_{\on{sc}})$} because
\begin{equation}\label{tildequotient}
\widetilde{\mu}^{(x)}_p:=\frac{\widetilde{M}^{(x)}_p}{\widetilde{M}^{(x)}_{p-1}}=(\mu^{(x)}_{4p-3}\cdots\mu^{(x)}_{4p})^{1/4},\;\;\;p\in\NN_{>0},\hspace{30pt}\widetilde{\mu}^{(x)}_0:=1,
\end{equation}
see also \cite[$(2.6)$]{subaddlike}.

\begin{lemma}\label{tildemarix}
Let $\mathcal{M}=\{M^{(x)}: x\in\mathcal{I}\}$ be \hyperlink{Msc}{$(\mathcal{M}_{\on{sc}})$}. Then we get
$$\forall\;x\in\mathcal{I}:\;\;\;M^{(x)}\le\widetilde{M}^{(x)},$$
and the following conditions are equivalent:
\begin{itemize}
\item[$(i)$] $\mathcal{M}$ is satisfying $\hyperlink{R-mg}{(\mathcal{M}_{\{\on{mg}\}})}$ resp. $\hyperlink{B-mg}{(\mathcal{M}_{(\on{mg})})}$,

\item[$(ii)$] $\widetilde{\mathcal{M}}$ and $\mathcal{M}$ are $R$- resp. $B$-equivalent.
\end{itemize}
\end{lemma}

\demo{Proof}
This follows immediately by a word-for-word repetition of the proof given in \cite[Lemma 2.2]{subaddlike} (with $C=4$ there) and by taking into account the equivalences between $(I_R)$ and $(II_R)$ resp. between $(I_B)$ and $(II_B)$.

Recall that for any $M\in\hyperlink{LCset}{\mathcal{LC}}$ and $C\in\NN_{>0}$ we get that
\begin{equation*}\label{rootgeneralizedincr}
\forall\;p\in\NN:\;\;\;(M_{Cp})^{1/C}\le(M_{(C+1)p})^{1/(C+1)},
\end{equation*}
because $p\mapsto(M_p)^{1/p}$ is nondecreasing by log-convexity and normalization.
\qed\enddemo

The importance of this auxiliary matrix $\widetilde{\mathcal{M}}$ is given by the following result.

\begin{lemma}\label{tildemariximplies}
Let $\mathcal{M}=\{M^{(x)}: x\in\mathcal{I}\}$ be \hyperlink{Msc}{$(\mathcal{M}_{\on{sc}})$} and $\widetilde{\mathcal{M}}$ shall denote the matrix defined in \eqref{tildedef}. Then we get:
\begin{equation}\label{tildemariximpliesequ}
\forall\;x\in\mathcal{I}\;\exists\;A\ge 1\;\forall\;p\in\NN:\;\;\;\mu^{(x)}_{2p}\le A\widetilde{\mu}^{(x)}_p.
\end{equation}
\end{lemma}

\demo{Proof}
By \eqref{tildequotient} we get $\widetilde{\mu}^{(x)}_p=(\mu^{(x)}_{4p-3}\cdots\mu^{(x)}_{4p})^{1/4}$ for all $p\ge 1$ and so
\begin{equation}\label{tildemariximpliesequ1}
\forall\;p\ge 2:\;\;\;\mu^{(x)}_{2p}\le\widetilde{\mu}^{(x)}_p\Leftrightarrow(\mu^{(x)}_{2p})^4\le\mu^{(x)}_{4p-3}\cdots\mu^{(x)}_{4p},
\end{equation}
which is valid because $p\mapsto\mu^{(x)}_p$ is nondecreasing by log-convexity and since $2p\le 4p-3\Leftrightarrow\frac{3}{2}\le p$. Hence, by choosing $A\ge 1$ sufficiently large, we have shown \eqref{tildemariximpliesequ}.
\qed\enddemo

Let $\mathcal{M}=\{M^{(x)}: x\in\mathcal{I}\}$ be \hyperlink{Msc}{$(\mathcal{M}_{\on{sc}})$}. Then we put (as a set) $$\mathcal{M}^{\on{mg}}:=\mathcal{M}\cup\widetilde{\mathcal{M}}.$$
Thus $\mathcal{M}\subseteq\mathcal{M}^{\on{mg}}$ and by Lemma \ref{tildemarix} we see:
\begin{itemize}
\item[$(*)$] $\mathcal{M}$ and $\mathcal{M}^{\on{mg}}$ are $R$- resp. $B$-equivalent if and only if $\mathcal{M}$ satisfies $\hyperlink{R-mg}{(\mathcal{M}_{\{\on{mg}\}})}$ resp. $\hyperlink{B-mg}{(\mathcal{M}_{(\on{mg})})}$.

\item[$(*)$] In this case we get $\mathcal{E}_{\{\mathcal{M}\}}=\mathcal{E}_{\{\widetilde{\mathcal{M}}\}}=\mathcal{E}_{\{\mathcal{M}^{\on{mg}}\}}$ resp. $\mathcal{E}_{(\mathcal{M})}=\mathcal{E}_{(\widetilde{\mathcal{M}})}=\mathcal{E}_{(\mathcal{M}^{\on{mg}})}$.

\item[$(*)$] However, note that $\mathcal{M}^{\on{mg}}$ is formally not a weight matrix as defined in Section \ref{classesweightmatrices} since the pointwise order $\le$ may fail in general.
\end{itemize}

\begin{remark}\label{tildemariximpliesrem}
Consider $\mathcal{M}\equiv\mathcal{M}_{\omega}$, with $\mathcal{M}_{\omega}$ being the matrix associated with a given weight $\omega\in\hyperlink{omset0}{\mathcal{W}_0}$. Then by Lemma \ref{tildemarix} and \eqref{newmoderategrowth} the matrices $\mathcal{M}_{\omega}$ and $\widetilde{\mathcal{M}_{\omega}}$ are both $R$- and $B$-equivalent.

More precisely, by definition we see
\begin{equation}\label{tildemariximpliesequforom}
\forall\;x>0\;\forall\;p\in\NN:\;\;\;\widetilde{W}^{(x)}_p=\left(W^{(x)}_{4p}\right)^{1/4}=W^{(4x)}_p.
\end{equation}
Consequently, \eqref{tildemariximpliesequ1} is consistent with \cite[Lemma 2.6]{whitneyextensionweightmatrix} and as sets one has $$\widetilde{\mathcal{M}_{\omega}}\subsetneq\mathcal{M}_{\omega}\;\;\;\text{and}\;\;\;\mathcal{M}_{\omega}^{\on{mg}}=\mathcal{M}_{\omega}.$$
However, by the $R$- and $B$-equivalence between $\widetilde{\mathcal{M}_{\omega}}$ and $\mathcal{M}_{\omega}$ we see that all matrices generate the same associated ultradifferentiable class of functions (as locally convex vector spaces), i.e.
$$\mathcal{E}_{\{\mathcal{M}_{\omega}\}}=\mathcal{E}_{\{\widetilde{\mathcal{M}_{\omega}}\}}=\mathcal{E}_{\{\mathcal{M}_{\omega}^{\on{mg}}\}}\;\;\;\text{and}\;\;\;\mathcal{E}_{(\mathcal{M}_{\omega})}=\mathcal{E}_{(\widetilde{\mathcal{M}_{\omega}})}=\mathcal{E}_{(\mathcal{M}_{\omega}^{\on{mg}})},$$
see \cite[Prop. 4.6]{compositionpaper}. If $\omega$ that has in addition \hyperlink{om1}{$(\omega_1)$}, then $\mathcal{E}_{\{\mathcal{M}_{\omega}\}}=\mathcal{E}_{\{\omega\}}$ and $\mathcal{E}_{(\mathcal{M}_{\omega})}=\mathcal{E}_{(\omega)}$ as well.
\end{remark}

By involving Lemma \ref{tildemariximplies} and the matrix $\mathcal{M}^{\on{mg}}$ we can connect the conditions listed in Propositions \ref{firstmg} and \ref{secondmg}.

\begin{theorem}\label{matrixmgtheorem}
Let $\mathcal{M}=\{M^{(x)}: x\in\mathcal{I}\}$ be \hyperlink{Msc}{$(\mathcal{M}_{\on{sc}})$}. Then the following are equivalent:
\begin{itemize}
\item[$(i)$] The matrix $\mathcal{M}$ satisfies $\hyperlink{R-mg}{(\mathcal{M}_{\{\on{mg}\}})}$ resp. $\hyperlink{B-mg}{(\mathcal{M}_{(\on{mg})})}$.

\item[$(ii)$] The matrix $\mathcal{M}^{\on{mg}}$ is $R$- resp. $B$-equivalent to $\mathcal{M}$ and all assertions $(I_R)-(V_R)$ resp. $(I_B)-(V_B)$ are satisfied for $\mathcal{M}^{\on{mg}}$.
\end{itemize}

Note that assertion $(i)$ is preserved under $R$- resp. $B$-equivalence of weight matrices.
\end{theorem}

\demo{Proof}
$(i)\Rightarrow(ii)$ By Proposition \ref{firstmg} the matrix $\mathcal{M}$ has $(I_R)-(III_R)$ resp. $(I_B)-(III_B)$ and, as pointed out before, by Lemma \ref{tildemarix} the matrices $\mathcal{M}$ and $\mathcal{M}^{\on{mg}}$ are $R$- resp. $B$-equivalent. Hence, also $\mathcal{M}^{\on{mg}}$ satisfies $(I_R)-(III_R)$ resp. $(I_B)-(III_B)$.

Moreover, note that \eqref{tildemariximpliesequ} precisely yields both $(IV_R)$ and $(IV_B)$ considered as a mixed condition between the matrices $\mathcal{M}$ and $\widetilde{\mathcal{M}}$, i.e. $(IV_R)$ and $(IV_B)$ holds for the matrix $\mathcal{M}^{\on{mg}}$.

Thus $(ii)$ is shown by combining Proposition \ref{firstmg} and \ref{secondmg}.\vspace{6pt}

$(ii)\Rightarrow(i)$ This follows immediately by the $R$- resp. $B$-equivalence and Lemma \ref{tildemarix} as commented before.
\qed\enddemo

In the weight function setting, by taking into account Remark \ref{tildemariximpliesrem}, we see that the technical matrices $\widetilde{\mathcal{M}_{\omega}}$ and $\mathcal{M}_{\omega}^{\on{mg}}$ are becoming superfluous and Theorem \ref{matrixmgtheorem} turns into the following ''more closed'' statement.

\begin{corollary}\label{matrixmgtheoremcor}
Let $\omega\in\hyperlink{omset0}{\mathcal{W}_0}$ be given and let $\mathcal{M}_{\omega}$ be the associated weight matrix. Then all assertions listed in Propositions \ref{firstmg} and \ref{secondmg} of both types are satisfied for $\mathcal{M}_{\omega}$.
\end{corollary}

\section{On the failure of the main characterizing result in the weight matrix setting}\label{failuresection}
The aim of this section is to see that for a given standard log-convex weight matrix the full generalization of Theorem \ref{mgreformulated} will fail. In particular, we prove that the generalization of \eqref{mgstrange} to the matrix setting is violated and we study this behavior in detail for the matrix $\mathcal{M}_{\omega_N}$, $N\in\hyperlink{LCset}{\mathcal{LC}}$.

\subsection{Quotient/root comparison conditions for weight matrices}
First let us introduce the matrix-type generalizations of \eqref{mgstrange}. We say that an abstractly given weight matrix $\mathcal{M}=\{M^{(x)}: x\in\mathcal{I}\}$ has the {\itshape quotient/root comparison property of Roumieu type}, if
\begin{equation}\label{rstrange}
\forall\;x\in\mathcal{I}\;\exists\;y\in\mathcal{I}\;\exists\;A\ge 1\;\forall\;p\in\NN_{>0}:\;\;\;\mu^{(x)}_p\le A(M^{(y)}_p)^{1/p},
\end{equation}
and of Beurling type, if
\begin{equation}\label{bstrange}
\forall\;x\in\mathcal{I}\;\exists\;y\in\mathcal{I}\;\exists\;A\ge 1\;\forall\;p\in\NN_{>0}:\;\;\;\mu^{(y)}_p\le A(M^{(x)}_p)^{1/p}.
\end{equation}

When $\mathcal{M}$ is \hyperlink{Msc}{$(\mathcal{M}_{\on{sc}})$}, then we get:

\begin{itemize}
\item[$(*)$] We can choose $x=y$ in \eqref{rstrange} and/or in \eqref{bstrange} if and only if $M^{(x)}$ has \hyperlink{mg}{$(\on{mg})$}, see Theorem \ref{mgreformulated}.

\item[$(*)$] When $\mathcal{M}=\{M\}$, $M\in\hyperlink{LCset}{\mathcal{LC}}$, then $\mathcal{M}$ has \eqref{rstrange} and/or \eqref{bstrange} if and only if $M$ has \hyperlink{mg}{$(\on{mg})$}.

\item[$(*)$] Similarly let us show now the following: When $\mathcal{M}$ is constant, then $\mathcal{M}$ has \eqref{rstrange} and/or \eqref{bstrange} if and only if each $M^{(x)}$ has \hyperlink{mg}{$(\on{mg})$}.

    On the one hand, if each $M^{(x)}$ has \hyperlink{mg}{$(\on{mg})$}, then the conclusion is clear by Theorem \ref{mgreformulated}.

    Conversely, since all sequences are equivalent, \eqref{rstrange} yields $\mu^{(x)}_p\le A(M^{(y)}_p)^{1/p}\le AB(M^{(x)}_p)^{1/p}$, i.e. \eqref{mgstrange} for $M^{(x)}$ and so \hyperlink{mg}{$(\on{mg})$} is verified.

    \eqref{bstrange} yields $\mu^{(y)}_p\le A(M^{(x)}_p)^{1/p}\le AB(M^{(y)}_p)^{1/p}$ and so \hyperlink{mg}{$(\on{mg})$} is verified for $M^{(y)}$. Since $M^{(y)}\hyperlink{approx}{\approx}M^{(x)}$ we get \hyperlink{mg}{$(\on{mg})$} for $M^{(x)}$, too.
\end{itemize}

\begin{remark}\label{strangeremark}
For $\mathcal{M}$ being \hyperlink{Msc}{$(\mathcal{M}_{\on{sc}})$} we gather some more information concerning these conditions:

\begin{itemize}
\item[$(i)$] In the Roumieu case, by \eqref{rstrange} we get that
    \begin{equation}\label{strangeremarkequ}
    \forall\;x\in\mathcal{I}\;\exists\;y\in\mathcal{I}\;\exists\;A\ge 1\;\forall\;p\in\NN_{>0}\;\forall\;z\ge y:\;\;\;\mu^{(x)}_p\le A(M^{(y)}_p)^{1/p}\le A(M^{(z)}_p)^{1/p}\le A\mu^{(z)}_p.
    \end{equation}
    \eqref{strangeremarkequ} yields that in the Roumieu case in \eqref{rstrange} w.l.o.g. we have $x\le y$ and we can restrict to $x,y\in\NN_{>0}$.

\item[$(ii)$] In the Beurling case, \eqref{bstrange} yields
\begin{equation}\label{strangeremarkequbeur}
    \forall\;x\in\mathcal{I}\;\exists\;y\in\mathcal{I}\;\exists\;A\ge 1\;\forall\;p\in\NN_{>0}\;\forall\;x'\ge x:\;\;\;\mu^{(y)}_p\le A(M^{(x)}_p)^{1/p}\le A(M^{(x')}_p)^{1/p}\le A\mu^{(x')}_p.
    \end{equation}
Here we are interested in small indices and so w.l.o.g. $x=\frac{1}{n}$, $n\in\NN_{>0}$. The choice $y>x$ in \eqref{bstrange} would give $(M^{(y)}_p)^{1/p}\le\mu^{(y)}_p\le A(M^{(x)}_p)^{1/p}$, hence $M^{(y)}\hyperlink{preceq}{\preceq} M^{(x)}\le M^{(y)}$ and so $M^{(y)}\hyperlink{approx}{\approx}M^{(x)}$. But, as $x\rightarrow 0$, this is only possible if $\mathcal{M}$ is constant or if the Beurling assertion in Remark \ref{constantremark} applies, i.e. the matrix stabilizes at some sufficiently small index. In the first case this means that each $M^{(x)}$ has \hyperlink{mg}{$(\on{mg})$} (as shown before) and in the second that all sequences with sufficiently small indices eventually have \hyperlink{mg}{$(\on{mg})$}.

Summarizing, in \eqref{bstrange} we can assume $y\le x$ and restrict to $y=\frac{1}{n_1}$ for some $n_1\in\NN_{>0}$, $n_1\ge n$.

\item[$(iii)$] By Stirling's formula we see that \eqref{rstrange} and \eqref{bstrange} are preserved under the mapping $\pi^s$, see \eqref{pitrafo}, i.e. \eqref{rstrange} and/or \eqref{bstrange} hold true simultaneously for some/any matrix $\{(p!^sM^{(x)}_p)_{p\in\NN}: x>0\}$, $s\in\RR$. However, for $s<0$ in general the resulting matrix will be not \hyperlink{Msc}{$(\mathcal{M}_{\on{sc}})$} anymore.
\end{itemize}
\end{remark}

In general \eqref{rstrange} resp. \eqref{bstrange} might be not preserved under $R$- resp. $B$-equivalence. However, when we consider for given matrices $\mathcal{M}=\{M^{(x)}: x>0\}$ and $\mathcal{N}=\{N^{(y)}: y>0\}$ the stronger relations
$$\forall\;x\in\mathcal{I}\;\exists\;y\in\mathcal{I}\;\exists\;A\ge 1\;\forall\;p\in\NN:\;\;\;\mu^{(x)}_p\le A\nu^{(y)}_p,$$
resp.
$$\forall\;x\in\mathcal{I}\;\exists\;y\in\mathcal{I}\;\exists\;A\ge 1\;\forall\;p\in\NN:\;\;\;\mu^{(y)}_p\le A\nu^{(x)}_p,$$
and if we denote the corresponding equivalence relation by $\{\cong\}$ resp. $(\cong)$, then it is immediate to see that \eqref{rstrange} is preserved under $\{\cong\}$ and \eqref{bstrange} under $(\cong)$: For the Roumieu case we get
\begin{align*}
&\forall\;x\in\mathcal{I}\;\exists\;y\in\mathcal{I}\;\exists\;z\in\mathcal{I}\;\exists\;z_1\in\mathcal{I}\;\exists\;A,B,C\ge 1\;\forall\;p\in\NN_{>0}:
\\&
\nu^{(x)}_p\le B\mu^{(y)}_p\le AB(M^{(z)}_p)^{1/p}\le ABC(N^{(z_1)}_p)^{1/p},
\end{align*}
and the Beurling case is similar.\vspace{6pt}

{\itshape Note:}

\begin{itemize}
\item[$(*)$] Even if we can choose $x=y$ in \eqref{rstrange} resp. \eqref{bstrange}, which amounts to the fact that each sequence has \hyperlink{mg}{$(\on{mg})$} (recall Theorem \ref{mgreformulated}), then in the general (non-constant) case it is not clear that automatically a weight matrix $\mathcal{N}$ which is $R$- resp. $B$-equivalent to $\mathcal{M}$ has \eqref{rstrange} resp. \eqref{bstrange}. (In general, $\mathcal{N}$ will have \hyperlink{R-mg}{$(\mathcal{M}_{\{\on{mg}\}})$} resp. \hyperlink{B-mg}{$(\mathcal{M}_{(\on{mg})})$}.)

\item[$(*)$] If $\mathcal{M}=\{M\}$, then there is no difference between $\{\cong\}$ resp. $(\cong)$ and $R$- resp. $B$-equivalence, i.e. equivalence: \eqref{rstrange} and \eqref{bstrange} are then precisely \eqref{mgstrange} for $M\in\hyperlink{LCset}{\mathcal{LC}}$ and so equivalent to \hyperlink{mg}{$(\on{mg})$} which is preserved under \hyperlink{approx}{$\approx$}. So each $N\in\hyperlink{LCset}{\mathcal{LC}}$ satisfying $N\hyperlink{approx}{\approx}M$ has \hyperlink{mg}{$(\on{mg})$} and so \eqref{mgstrange} too (see also \cite[Prop. 2.7, Rem. 2.8]{JimenezGarridoSanz}).
\end{itemize}

Next we show that by using these conditions we can prove Theorem \ref{matrixmgtheorem} without involving the technical auxiliary matrices $\widetilde{\mathcal{M}}$ and $\mathcal{M}^{\on{mg}}$.

\begin{proposition}\label{thirdmg}
Let $\mathcal{M}=\{M^{(x)}: x\in\mathcal{I}\}$ be \hyperlink{Msc}{$(\mathcal{M}_{\on{sc}})$}. Then we get:
\begin{itemize}
\item[$(i)$] If $\mathcal{M}$ has in addition \eqref{rstrange}, then all Roumieu type assertions $(I_R)-(V_R)$ in Propositions \ref{firstmg} and \ref{secondmg} are equivalent.

\item[$(ii)$] If $\mathcal{M}$ has in addition \eqref{bstrange}, then all Beurling type assertions $(I_B)-(V_B)$ in Propositions \ref{firstmg} and \ref{secondmg} are equivalent.
\end{itemize}
\end{proposition}

Note that this result is slightly stronger than Theorem \ref{matrixmgtheorem} (and Corollary \ref{matrixmgtheoremcor}) since there we get that $(I_R)-(V_R)$ and/or $(I_B)-(V_B)$ are satisfied simultaneously, which is sufficient knowledge for applications, however the equivalence of all five conditions of the particular type has not been shown since the implications $(II_R)\Rightarrow(IV_R)$ and $(II_B)\Rightarrow(IV_B)$ have not been verified for $\mathcal{M}$ directly.

\demo{Proof}
We only consider the Roumieu case in detail: When $\mathcal{M}$ is satisfying $(II_R)$, then by combining this with \eqref{rstrange} and the fact that each sequence is normalized we get:
$$\forall\;x\in\mathcal{I}\;\exists\;y,z\in\mathcal{I}\;\exists\;A,C\ge 1\;\forall\;p\in\NN_{>0}:\;\;\;\mu^{(x)}_{2p}\le A(M^{(y)}_{2p})^{1/(2p)}\le AC(M^{(z)}_p)^{1/p}\le AC\mu^{(z)}_p.$$
Hence $(IV_R)$ is verified (the case $p=0$ yields $\mu^{(x)}_0=1=\mu^{(y)}_0$ for any indices $x,y\in\mathcal{I}$).
\qed\enddemo

\begin{remark}
Let $\mathcal{M}=\{M^{(x)}: x\in\mathcal{I}\}$ be \hyperlink{Msc}{$(\mathcal{M}_{\on{sc}})$}. By inspecting the proof of \cite[Appendix B, $(2)\Leftrightarrow(4)$]{matsumotopseudo} and \cite[Lemma 5.3 $(4)\Rightarrow(3)$]{PetzscheVogt} we see that the quotient/root comparison properties of the particular type is related to $(II_R)$ resp. $(II_B)$ as follows:

\begin{itemize}
\item[$(i)$] If $\mathcal{M}$ has \eqref{rstrange}, then by applying log-convexity for $M^{(x)}$ we can estimate by
$$\frac{M^{(x)}_{2p}}{M^{(x)}_p}=\mu^{(x)}_{p+1}\cdots\mu^{(x)}_{2p}\le(\mu^{(x)}_{2p})^{p}\le A^{p}(M^{(y)}_{2p})^{1/2},$$
consequently $\mathcal{M}$ satisfies
$$\forall\;x\in\mathcal{I}\;\exists\;y\in\mathcal{I}\;\exists\;A\ge 1\;\forall\;p\in\NN:\;\;\;\frac{M^{(x)}_{2p}}{(M^{(y)}_{2p})^{1/2}}\le A^pM^{(x)}_p.$$

\item[$(ii)$] If $\mathcal{M}$ has $(II_R)$, then
$$(\mu^{(x)}_p)^p\le\mu^{(x)}_{p+1}\cdots\mu^{(x)}_{2p}=\frac{M^{(x)}_{2p}}{M^{(x)}_p}\le A^{2p}\frac{(M^{(y)}_p)^2}{M^{(x)}_p},$$
consequently $\mathcal{M}$ satisfies
$$\forall\;x\in\mathcal{I}\;\exists\;y\in\mathcal{I}\;\exists\;A\ge 1\;\forall\;p\in\NN_{>0}:\;\;\;\mu^{(x)}_p\le A\left(\frac{(M^{(y)}_p)^2}{M^{(x)}_p}\right)^{1/p}.$$
\end{itemize}

The statements for the Beurling case are analogous.
\end{remark}

\subsection{Weight matrices associated with associated weight functions}
We are interested now in studying the situation when the matrix is associated with a weight function $\omega$. Recall that in this case by Corollary \ref{matrixmgtheoremcor} all assertions listed in Propositions \ref{firstmg} and \ref{secondmg} hold true for $\mathcal{M}_{\omega}$ and so, if $\mathcal{M}_{\omega}$ has in addition \eqref{rstrange} resp. \eqref{bstrange} then we can apply Proposition \ref{thirdmg}. Recall that by the convexity and increasing properties of $\varphi^{*}_{\omega}$ we even have $\vartheta^{(x)}\le\vartheta^{(y)}$ for $x\le y$, see \cite[Sect. 2.5]{whitneyextensionweightmatrix}.

As a special but interesting and concrete case we focus now on $\omega\equiv\omega_M$, $M\in\hyperlink{LCset}{\mathcal{LC}}$.\vspace{6pt}

First we characterize \eqref{rstrange} and \eqref{bstrange} for $\mathcal{M}_{\omega_M}$ in terms of a growth property for $M$. More precisely we prove the following result.

\begin{proposition}\label{strangechar}
Let $M\in\hyperlink{LCset}{\mathcal{LC}}$ be given and let $\mathcal{M}_{\omega_M}=\{M^{(l)}: l>0\}$ be the matrix associated with $\omega_M$. Then the following are equivalent:
\begin{itemize}
\item[$(a_R)$] $\mathcal{M}_{\omega_M}$ satisfies
\begin{equation}\label{rstrangemod}
\exists\;A\ge 1\;\exists\;c\ge 1\;\forall\;x\in\NN_{>0}\;\forall\;p\in\NN_{>0}:\;\;\;\mu^{(x)}_p\le A(M^{(cx)}_p)^{1/p}.
\end{equation}

\item[$(b_R)$] $\mathcal{M}_{\omega_M}$ has the {\itshape quotient/root comparison property of Roumieu type} \eqref{rstrange}.

\item[$(a_B)$] $\mathcal{M}_{\omega_M}$ satisfies
\begin{equation}\label{bstrangemod}
\exists\;c\ge 1\;\forall\;x\in\NN_{>0}\;\exists\;A\ge 1\;\;\forall\;p\in\NN_{>0}:\;\;\;\mu^{(1/(cx))}_p\le A(M^{(1/x)}_p)^{1/p}.
\end{equation}

\item[$(b_B)$] $\mathcal{M}_{\omega_M}$ has the {\itshape quotient/root comparison property of Beurling type} \eqref{bstrange}.

\item[$(c)$] $M$ satisfies
\begin{equation}\label{genmg}
\exists\;d\in\NN_{>0}\;\exists\;A\ge 1\;\forall\;p\in\NN_{>0}:\;\;\;\mu_p\le A(M_{dp})^{1/(dp)}.
\end{equation}
\end{itemize}
\end{proposition}

{\itshape Conclusion:} By combining this result with Corollary \ref{matrixmgtheoremcor} we get that the matrix $\mathcal{M}_{\omega_M}$ satisfies all matrix-type-generalizations of the assertions listed in Theorem \ref{mgreformulated} (for both types) simultaneously if and only if $M$ satisfies \eqref{genmg}.

Note that:

\begin{itemize}
\item[$(i)$] The mapping $d\mapsto(M_{dp})^{1/(dp)}$ is nondecreasing for all $p\in\NN_{>0}$ arbitrary but fixed. So, if \eqref{genmg} is valid for some $d_0$, then for all $d\ge d_0$ as well.

\item[$(ii)$] By Stirling's formula, \eqref{genmg} is preserved under the mapping $\pi^s$, $s\in\RR$, see \eqref{pitrafo}, by only changing the constant $A$ but with the same choice for $d$. However, when $s<0$ then in general one might lose log-convexity for the resulting sequence and also $(p!^sM_p)^{1/p}\rightarrow+\infty$ as $p\rightarrow+\infty$ might fail.

\item[$(iii)$] By $(i)$ and $(ii)$ in Remark \ref{strangeremark}, in \eqref{rstrange} we can restrict ourselves w.l.o.g. to $x,y\in\NN_{>0}$ with $x\le y$, whereas in \eqref{bstrange} to $x=1/n$ and $y=1/n_1$ with $n,n_1\in\NN_{>0}$ and $n_1\ge n$.
\end{itemize}

\demo{Proof}
$(a_R)\Rightarrow(b_R)$ This is clear since in \eqref{rstrangemod} both $A$ and $c$ are not depending on given index $x$ and by taking into account comment $(iii)$ before.\vspace{6pt}

$(b_R)\Rightarrow(c)$ By assumption $\mathcal{M}_{\omega_M}$  has \eqref{rstrange} and by choosing $x=1$ there we get
$$\exists\;y\in\NN_{>0}\;\exists\;A\ge 1\;\forall\;p\in\NN_{>0}:\;\;\;\mu^{(1)}_p=\mu_p\le A(M^{(y)}_p)^{1/p}=A(M^{(1)}_{yp})^{1/(yp)}=A(M_{yp})^{1/(yp)},$$
where we have used \eqref{Mgoodtransform} and \eqref{Mxonestable}. Thus \eqref{genmg} is shown with $A=A$ and $d:=y$.\vspace{6pt}

$(c)\Rightarrow(a_R)$ First we show the following:
\begin{equation}\label{beta1forx}
\forall\;x\in\NN_{>0}\;\forall\;p\in\NN:\;\;\;\mu^{(x)}_p\le\mu_{xp}.
\end{equation}
If $p=0$, then we have the equality $\mu^{(x)}_0=1=\mu_0$. For all $p,x\in\NN_{>0}$, by using again \eqref{Mgoodtransform} and \eqref{Mxonestable}, we get that
\begin{equation}\label{beta1forx1}
\mu^{(x)}_p:=\frac{M^{(x)}_p}{M^{(x)}_{p-1}}=\frac{(M^{(1)}_{xp})^{1/x}}{(M^{(1)}_{xp-x})^{1/x}}=\left(\frac{M_{xp}}{M_{xp-x}}\right)^{1/x}=(\mu_{xp-x+1}\cdots\mu_{xp})^{1/x}\le(\mu_{xp})^{x/x}=\mu_{xp},
\end{equation}
where the last estimate holds by log-convexity for $M$. Combining \eqref{beta1forx} with \eqref{genmg} it follows that
$$\exists\;A\ge 1\;\exists\;d\in\NN_{>0}\;\forall\;x,p\in\NN_{>0}:\;\;\;\mu^{(x)}_p\le\mu_{xp}\le A(M_{dxp})^{1/(dxp)}=A(M^{(1)}_{dxp})^{1/(dxp)}=A(M^{(dx)}_p)^{1/p},$$
hence \eqref{rstrangemod} is verified with $A=A$ and $c:=d$.

Note that both $A$ and $d$ are only depending on given sequence $M$ but not on the index $x$.\vspace{12pt}

$(a_B)\Rightarrow(b_B)$ This is clear since $c$ is not depending on $x$ (and by taking into account again comment $(iii)$).\vspace{6pt}

$(b_B)\Rightarrow(c)$ First we are showing the following formula which is ''dual'' to \eqref{beta1forx}:
\begin{equation}\label{beta1forxbeur1}
\forall\;y\in\NN_{>0}\;\forall\;p\in\NN:\;\;\;\mu^{(1/y)}_{yp}\ge\mu^{(1)}_p=\mu_p.
\end{equation}
For $y=1$ we trivially get equality for all $p\in\NN$ and also for $p=0$ because then $\mu^{(1/y)}_{0}=1=\mu_0$. So let $y\ge 2$ and $p\ge 1$, then by definition, the log-convexity for $M^{(1/y)}$ and again by \eqref{Mgoodtransform} and \eqref{Mxonestable} we obtain
\begin{align*}
\mu^{(1/y)}_{yp}&=\frac{M^{(1/y)}_{yp}}{M^{(1/y)}_{yp-1}}=\frac{M^{(1/y)}_{yp}}{M^{(1/y)}_{yp-y}}\frac{M^{(1/y)}_{yp-2}}{M^{(1/y)}_{yp-1}}\frac{M^{(1/y)}_{yp-3}}{M^{(1/y)}_{yp-2}}\cdots\frac{M^{(1/y)}_{yp-y}}{M^{(1/y)}_{yp-y+1}}
\\&
=\frac{M^{(1/y)}_{yp}}{M^{(1/y)}_{yp-y}}\frac{1}{\mu^{(1/y)}_{yp-1}}\frac{1}{\mu^{(1/y)}_{yp-2}}\cdots\frac{1}{\mu^{(1/y)}_{yp-y+1}}\ge\frac{M^{(1/y)}_{yp}}{M^{(1/y)}_{yp-y}}\left(\frac{1}{\mu^{(1/y)}_{yp}}\right)^{y-1}
\\&
=\left(\frac{M^{(1)}_{p}}{M^{(1)}_{p-1}}\right)^y\left(\frac{1}{\mu^{(1/y)}_{yp}}\right)^{y-1}=(\mu_p)^y\left(\frac{1}{\mu^{(1/y)}_{yp}}\right)^{y-1}.
\end{align*}

We choose $x=1$ in \eqref{bstrange} and get
$$\exists\;y\in\NN_{>0}\;\exists\;A\ge 1\;\forall\;p\in\NN_{>0}:\;\;\;(M_p)^{1/p}=(M^{(1)}_p)^{1/p}\ge\frac{1}{A}\mu^{(1/y)}_p,$$
hence by taking into account \eqref{beta1forxbeur1}
$$\exists\;y\in\NN_{>0}\;\exists\;A\ge 1\;\forall\;p\in\NN_{>0}:\;\;\;(M_{yp})^{1/(yp)}\ge\frac{1}{A}\mu^{(1/y)}_{yp}\ge\frac{1}{A}\mu_p.$$
Consequently, we have verified \eqref{genmg} with $d:=y$ and $A=A$.\vspace{6pt}

$(c)\Rightarrow(a_B)$ Similarly, as in \eqref{beta1forx} we get
\begin{equation}\label{beta1forxbeur}
\forall\;x\in\NN_{>0}\;\forall\;c\in\NN_{>0}\;\forall\;p\in\NN_{>0}:\;\;\;\mu^{(x)}_{p}\ge\mu^{(x/c)}_{c(p-1)},
\end{equation}
because
$$(\mu^{(x)}_{p})^{c}=\left(\frac{M^{(x)}_{p}}{M^{(x)}_{p-1}}\right)^c=\frac{M^{(x/c)}_{cp}}{M^{(x/c)}_{c(p-1)}}=\mu^{(x/c)}_{cp-c+1}\cdots\mu^{(x/c)}_{cp}\ge\left(\mu^{(x/c)}_{c(p-1)}\right)^c.$$
By applying \eqref{beta1forxbeur} to $x=1$ and by involving \eqref{genmg} we can find $A\ge 1$ and $d\in\NN_{>0}$ such that for any $x,c,p\in\NN_{>0}$ we get the following estimation:
\begin{align*}
\mu^{(1/c)}_{c(p-1)}\le\mu^{(1)}_p=\mu_p\le A(M_{dp})^{1/(dp)}=A(M^{(1)}_{dp})^{1/(dp)}=A(M^{(1/x)}_{xdp})^{1/(xdp)}.
\end{align*}
So let $x\in\NN_{>0}$ be arbitrary but from now on fixed. Then set $c=2xd$ and so for all $q\in\NN_{>0}$, $q\ge 2$, we have
\begin{align*}
(M^{(1/x)}_{xdq})^{1/(xdq)}\ge\frac{1}{A}\mu^{(1/(2xd))}_{2xd(q-1)}\ge\frac{1}{A}\mu^{(1/(2xd))}_{xdq},
\end{align*}
which holds by log-convexity for $M^{1/(2xd)}$ and since $2xd(q-1)\ge xdq\Leftrightarrow q\ge 2$. This proves \eqref{bstrangemod} for the choice $c:=2d$, $A=A$ and all $p\in\NN_{>0}$ such that $p=xdq$ for some $q\in\NN_{>0}$, $q\ge 2$. Next let $xdq<p<xd(q+1)$, $q\ge 2$, and since $p\mapsto(M^{(l)}_p)^{1/p}$ is nondecreasing for each $l>0$ we estimate as follows:
\begin{equation}\label{strangecharequ}
(M^{(1/x)}_p)^{1/p}\ge(M^{(1/x)}_{xdq})^{1/(xdq)}\ge\frac{1}{A}\mu^{(1/(2xd))}_{xdq}\ge\frac{1}{AB}\mu^{(1/(8xd))}_{2xdq}\ge\frac{1}{AB}\mu^{(1/(8xd))}_{xd(q+1)}\ge\frac{1}{AB}\mu^{(1/(8xd))}_{p}.
\end{equation}
Apart from the log-convexity for $M^{(1/(8xd))}$ and $2xdq\ge xd(q+1)\Leftrightarrow q\ge 1$ we have used \eqref{tildemariximpliesequ} and \eqref{tildemariximpliesequforom} and so the constant $B$ is here also depending on given index $x$. By increasing the constant $B$ (if necessary) we finally get \eqref{strangecharequ} also for $1\le p<2xd$. Thus \eqref{bstrangemod} is verified when choosing the constant $AB$ (depending on given $x$) and for $c$ we set $c:=8d$ not depending on the index $x$.
\qed\enddemo

\begin{remark}\label{genmgremark}
Let $M\in\hyperlink{LCset}{\mathcal{LC}}$, we give now several comments on the new characterizing condition \eqref{genmg}. Since $p\mapsto(M_p)^{1/p}$ is nondecreasing, \eqref{genmg} can be viewed as an {\itshape almost/generalized moderate growth condition} for given $M$.\vspace{6pt}

\begin{itemize}
\item[$(i)$] First, we recall the following equivalences (by combining Theorem \ref{mgreformulated}, Lemma \ref{assoweightomega0} and $(iii)$ in Section \ref{classesweightmatrices}):

\begin{itemize}
\item[$(*)$] \eqref{genmg} holds with $d=1$, i.e. \eqref{mgstrange} is valid,

\item[$(*)$] Theorem \ref{mgreformulated} applies to $M$,

\item[$(*)$] Theorem \ref{mgreformulated} applies to some/each sequence $M^{(x)}$,

\item[$(*)$] $M^{(x)}\hyperlink{approx}{\approx}M^{(y)}$ for all $x,y>0$, i.e. $\mathcal{M}_{\omega_M}$ is constant (in particular, each $M^{(x)}$ is equivalent to $M\equiv M^{(1)}$).
\end{itemize}

\item[$(ii)$] \eqref{genmg} for $d\ge 2$ can be satisfied for very fast increasing sequences (in particular when \hyperlink{mg}{$(\on{mg})$} is violated). For this let $q>1$ and consider $M_p:=q^{p^n}$ for $n\in\NN_{\ge 2}$. Thus for $n=2$ we obtain the so-called $q${\itshape -Gevrey sequences.} Then for $p\in\NN_{>0}$ one has $$\mu_p=q^{p^n-(p-1)^n}=q^{p^n-\sum_{k=0}^n\binom{n}{k}p^k(-1)^{n-k}}=q^{-\sum_{k=0}^{n-1}\binom{n}{k}p^k(-1)^{n-k}},$$
    and so $\mu_p=q^{np^{n-1}+O(p^{n-2})}$, whereas $(M_{dp})^{1/(dp)}=q^{(dp)^{n-1}}$. Hence \eqref{genmg} follows, e.g. for $n=2$ we get $\mu_p=q^{2p-1}\le(M_{dp})^{1/(dp)}=q^{dp}$ for all $p\in\NN$ by choosing some $d\ge 2$.\vspace{6pt}

Note that none of such sequences has \hyperlink{mg}{$(\on{mg})$}, for $n\ge 3$ even derivation closedness is violated since $\mu_p\le D^p$ for all $p\in\NN$ resp. $\log(M_p)=O(p^2)$ as $p\rightarrow\infty$ fails for all such values. Those sequences are arising in the weight matrices of the weights $\omega_s(t):=\max\{0,\log(t)^s\}$, $s>1$, see \cite[Sect. 5.5]{whitneyextensionweightmatrix} and the references therein.

\item[$(iii)$] But even for much faster growing sequences we can have \eqref{genmg} with $d=2$: Let
$$M_p:=e^{e^p},\;\;\;p\ge 1,\hspace{20pt}M_0:=1,$$
then $(M_{dp})^{1/(dp)}=e^{e^{dp}/(dp)}$ for all $d,p\ge 1$ and $\mu_p=e^{e^p-e^{p-1}}=e^{e^{p-1}(e-1)}$ for $p\ge 2$ and $\mu_1=M_1=e^e$. Hence for all $p\ge 2$
\begin{align*}
&\mu_p\le A(M_{dp})^{1/(dp)}\Leftrightarrow e^{p-1}(e-1)\le\log(A)+\frac{e^{dp}}{dp}
\\&
\Leftrightarrow dp(e^{-1}(e-1)-\log(A)/e^{p})=dp(1-e^{-1}-\log(A)/e^p)\le e^{p(d-1)},
\end{align*}
which holds true for all $p\ge 2$ even with $d=2$ by taking any $A\ge 1$ (even the choice $A=1$ is sufficient). For $p=1$ we have with $d=2$ that
$$e^e=\mu_1\le(M_2)^{1/2}=e^{e^2/2}.$$
\end{itemize}
\end{remark}

Summarizing, the new characterizing condition \eqref{genmg} motivates the following definition:\vspace{6pt}

For given $M\in\hyperlink{LCset}{\mathcal{LC}}$ (or even assume $M\in\RR_{>0}^{\NN}$) let the {\itshape moderate growth index} $g(M)$ be defined by
$$g(M):=\min\{d\in\NN_{>0}:\;\;\;\eqref{genmg}\;\text{holds true}\},$$
and let us set $g(M):=+\infty$ if \eqref{genmg} is violated. So

\begin{itemize}
\item[$(i)$] $g(M)=1$ holds if and only if $M$ has \hyperlink{mg}{$(\on{mg})$}.

\item[$(ii)$] $g(M)<+\infty$ if and only if $\mathcal{M}_{\omega_M}$ has \eqref{rstrange} and/or \eqref{bstrange}, i.e. if and only if the matrix $\mathcal{M}_{\omega_M}$ satisfies all generalizations of the assertions listed in Theorem \ref{mgreformulated} for both types.

\item[$(iii)$] Note that this index is preserved under the mapping $\pi^s$, $s\in\RR$.
\end{itemize}

Next we study how \eqref{genmg} is transformed under the equivalence relation between weight sequences.

\begin{lemma}\label{equlemma}
Let $M,N\in\RR_{>0}^{\NN}$ be given with $M_0=N_0=1$ and assume that $M\hyperlink{approx}{\approx}N$. Moreover assume that $M$ has \eqref{genmg}, then $N$ has to satisfy
\begin{equation}\label{equlemma1equ1}
\exists\;d\in\NN_{>0}\;\exists\;A,C\ge 1\;\forall\;p\in\NN_{>0}:\;\;\;\nu_p\le AC^{2p}(N_{dp})^{1/(dp)}\Leftrightarrow N_p\le AC^{2p}(N_{dp})^{1/(dp)}N_{p-1},
\end{equation}
where $d$ and $A$ are the parameters arising in \eqref{genmg}, i.e. only depending on given $M$, and $C$ is coming from the equivalence \hyperlink{approx}{$\approx$}, hence also depending on $N$.
\end{lemma}

\demo{Proof}
The equivalence between $M$ and $N$ yields
$$\exists\;C\ge 1\;\forall\;p\in\NN_{>0}:\;\;\;\mu_p=\frac{M_p}{M_{p-1}}\ge\frac{N_p}{C^p}\frac{1}{C^{p-1}N_{p-1}}=\frac{1}{C^{2p-1}}\nu_p.$$

Consequently we get $\frac{1}{C^{2p-1}}\nu_p\le\mu_p\le A(M_{dp})^{1/(dp)}\le AC(N_{dp})^{1/(dp)}$ and so \eqref{equlemma1equ1} is shown.
\qed\enddemo

Note that:

\begin{itemize}
\item[$(i)$] It is straightforward to see that \eqref{equlemma1equ1} is stable under equivalence of weight sequences and it is also preserved under the mapping $\pi^s$, $s\in\RR$ (uniformly with the same choice for $d$).

\item[$(ii)$] If $N\in\hyperlink{LCset}{\mathcal{LC}}$, then \hyperlink{dc}{$(\on{dc})$} does immediately imply \eqref{equlemma1equ1}.


\item[$(iii)$] However, the converse implication fails in general: For this consider the sequence in $(iii)$ in Remark \ref{genmgremark}, so \eqref{equlemma1equ1} follows with $d=2$ and $A=C=1$, but \hyperlink{dc}{$(\on{dc})$} fails.
\end{itemize}




\subsection{Weight matrices associated with general weight functions}
We turn now to the weight function situation and the matrix $\mathcal{M}_{\omega}$.

\begin{proposition}\label{strangecharom1}
Let $\omega\in\hyperlink{omset0}{\mathcal{W}_0}$ be given and let $\mathcal{M}_{\omega}=\{W^{(l)}: l>0\}$ be the associated weight matrix.

Then the following are equivalent:

\begin{itemize}
\item[$(i)$] $\mathcal{M}_{\omega}$ satisfies the {\itshape quotient/root comparison property of Roumieu type} \eqref{rstrange} and/or of {\itshape Beurling type} \eqref{bstrange},

\item[$(ii)$] we have
\begin{equation*}\label{genmgom1}
\exists\;A\ge 1\;\exists\;c\ge 1\;\forall\;x\in\NN_{>0}\;\forall\;p\in\NN_{>0}:\;\;\;\vartheta^{(x)}_p\le A(W^{(cx)}_p)^{1/p},
\end{equation*}

\item[$(iii)$] we have
\begin{equation*}
\exists\;c\ge 1\;\forall\;x\in\NN_{>0}\;\exists\;A\ge 1\;\;\forall\;p\in\NN_{>0}:\;\;\;\vartheta^{(1/(cx))}_p\le A(W^{(1/x)}_p)^{1/p},
\end{equation*}

\item[$(iv)$] we have
\begin{equation*}\label{genmgom1}
\exists\;d\in\NN_{>0}\;\exists\;A\ge 1\;\forall\;p\in\NN_{>0}:\;\;\;\vartheta^{(1)}_p\le A(W^{(1)}_{dp})^{1/(dp)},
\end{equation*}

\item[$(v)$] the matrix $\mathcal{M}_{\omega_{W^{(1)}}}$ associated with the weight $\omega_{W^{(1)}}\in\hyperlink{omset0}{\mathcal{W}_0}$ satisfies the {\itshape quotient/root comparison property of Roumieu type} \eqref{rstrange} and/or of {\itshape Beurling type} \eqref{bstrange}.
\end{itemize}
\end{proposition}

{\itshape Conclusion:} The equivalence between $(i)$ and $(v)$ yields that the information about the desired properties \eqref{rstrange} and/or \eqref{bstrange} for $\mathcal{M}_{\omega}$ is already encoded in a matrix associated with an associated weight function, i.e. in the single sequence $W^{(1)}$.

\demo{Proof}
$(i)\Leftrightarrow(ii)\Leftrightarrow(iii)\Leftrightarrow(iv)$ This follows by a word-for-word repetition of the proof given in Proposition \ref{strangechar} by using $M^{(x)}\equiv W^{(x)}$ and $\mu^{(x)}\equiv\vartheta^{(x)}$, see Section \ref{classesweightmatrices}. (We only have to skip the additional information $W^{(1)}\equiv M$, $\vartheta^{(1)}\equiv\mu$, see \eqref{Mxonestable}.)\vspace{6pt}

$(iv)\Leftrightarrow(v)$ This follows directly by applying Proposition \ref{strangechar} to $W^{(1)}\equiv M$.
\qed\enddemo

\subsection{A counter-example}\label{countersection}
We are now constructing a (counter)-example for which the matrix $\mathcal{M}_{\omega_M}$ associated with $\omega_M$ violates both \eqref{rstrange} and \eqref{bstrange}. So, in the general (non-constant) case, the generalization of Theorem \ref{mgreformulated} (\cite[Lemma 2.2]{whitneyextensionweightmatrix}) to the mixed setting fails.

\begin{theorem}\label{counter1}
There exist $N\in\hyperlink{LCset}{\mathcal{LC}}$ such that $N$ violates \eqref{genmg} (i.e. $g(N)=+\infty$).

By Proposition \ref{strangechar} this is equivalent to the fact that the matrix $\mathcal{M}_{\omega_N}$ associated with $\omega_N$ does not have \eqref{rstrange} and/or \eqref{bstrange}.

In particular, $N$ violates \hyperlink{mg}{$(\on{mg})$} equivalently $\mathcal{M}_{\omega_N}$ is non-constant, see $(i)$ in Remark \ref{genmgremark}.\vspace{6pt}

In addition, we can obtain some more properties:
\begin{itemize}
\item[$(a)$] $N$ can be chosen to be strong non-quasianalytic, i.e. $N$ is satisfying \hyperlink{beta1}{$(\beta_1)$}. In this case, $\omega_N$ has \eqref{assostrongnq} (strong non-quasianalyticity condition) as well, see \cite[Proposition 4.4]{Komatsu73}, which means that $\omega_N$ is even a {\itshape strong weight} in the notion of \cite{BonetBraunMeiseTaylorWhitneyextension}.

On the other hand $N$ can also be chosen to be quasianalytic, i.e. $N$ is satisfying
$$\sum_{p\ge 1}\frac{1}{\nu_p}=+\infty.$$

\item[$(b)$] $N$ can be chosen such that even \eqref{equlemma1equ1} is violated.
\end{itemize}
\end{theorem}

\demo{Proof}
We start with the following construction. Let $N$ be given by
\begin{equation}\label{Ndefidea}
N_p:=\exp(f(p)),\hspace{20pt}f:[0,+\infty)\rightarrow[0,+\infty),
\end{equation}
with $f$ being the convex (continuous) function defined as follows: The graph of $f$ is consisting of all straight lines connecting the points $\{(a_j,f(a_j)): j\ge 1\}$ with $(a_j)_{j\ge 1}$ being a sequence in $\NN$. The slope of the (straight) line connecting the points $(a_j,f(a_j))$ and $(a_{j+1},f(a_{j+1}))$ is denoted by $b_j\in\RR_{>0}$. By using this construction in order to get $N\in\hyperlink{LCset}{\mathcal{LC}}$ we require that $j\mapsto b_j$ is (strictly) increasing and such that $\lim_{j\rightarrow\infty}b_j=+\infty$. The strategy is to construct this graph step-wise subject to growth restrictions on $(a_j)_j$ and $(b_j)_j$. So, when given the values $a_j$, $f(a_j)$, $a_{j+1}$ and $b_j$, then we set $f(a_{j+1}):=b_j(a_{j+1}-a_j)+f(a_j)$ which determines automatically $f(p)$ for $a_j+1\le p\le a_{j+1}-1$.\vspace{6pt}

For the sequence $(a_j)_{j\ge 1}$ let us assume
\begin{equation}\label{aproperty}
a_1:=0,\;\;\;\;j(a_j+1)\le a_{j+1},\;\;\;\forall\;j\ge 1.
\end{equation}
First we set
$$f(a_1)=f(0):=0,$$
hence $N_0=1$ follows, and then
\begin{equation}\label{bproperty}
b_1>0\;\;\text{arbitrary},\;\;\;\;b_j>b_{j-1},\;\;\;\;b_j\ge\frac{j^2(a_j+1)+f(a_j)}{a_j},\;\;\;j\ge 2.
\end{equation}
This choice for $j\ge 2$ is possible since the value $f(a_j)$ is only depending on given points $a_1,\dots ,a_j$ and slopes $b_1,\dots,b_{j-1}$ and so $\lim_{j\rightarrow\infty}b_j=+\infty$ follows, too. Thus $N\in\hyperlink{LCset}{\mathcal{LC}}$ is verified.\vspace{6pt}

$N$ satisfies \eqref{genmg} if and only if
$$\exists\;d\in\NN_{>0}\;\exists\;A\ge 1\;\forall\;p\in\NN_{>0}:\;\;\;dp(f(p)-f(p-1)-\log(A))\le f(dp).$$
We prove that this condition is impossible.

Let now $d\in\NN_{>0}$, $d\ge 2$, and $A\ge 1$ be given, arbitrary (large) but from now on fixed. Then for all $l\ge d\ge 2$ the choice $p=a_l+1$ yields $f(p)-f(p-1)=b_l$ and
$$f(dp)=f(a_l)+b_l(dp-a_l)=f(a_l)+b_l(d(a_l+1)-a_l),$$
because $d(a_l+1)\le l(a_l+1)\le a_{l+1}$ by the choice of $(a_l)_l$, see \eqref{aproperty}.

Thus
\begin{align*}
&d(a_l+1)(b_l-\log(A))\le f(a_l)+b_l(d(a_l+1)-a_l)\Leftrightarrow d(a_l+1)\log(A)\ge-f(a_l)+a_lb_l
\\&
\Leftrightarrow\frac{d(a_l+1)\log(A)+f(a_l)}{a_l}\ge b_l,
\end{align*}
a contradiction as $l\rightarrow\infty$ by the choices of the numbers $b_l$, see \eqref{bproperty}.
\vspace{12pt}

However, it is not clear that $N$ has $(a)$ and/or $(b)$.

In order to get the strong non-quasianalyticity part of $(a)$ we recall the following technique, see also \cite[Example 7.10]{sectorialextensions} where a similar construction for $N$ has been considered:\vspace{6pt}

Since $N$ is log-convex we clearly get $\liminf_{p\rightarrow\infty}\frac{\nu_{Qp}}{\nu_p}\ge 1$ for any $Q\in\NN_{>0}$. But note that \hyperlink{beta1}{$(\beta_1)$}, i.e. $\liminf_{p\rightarrow\infty}\frac{\nu_{Qp}}{\nu_p}>Q$  and even $\liminf_{p\rightarrow\infty}\frac{\nu_{Qp}}{\nu_p}>1$ is violated for any $Q\in\NN_{\ge 2}$ by the requirement $j(a_j+1)\le a_{j+1}$ for all $j\ge 1$: We have $\nu_p=\exp(b_i)$ for all $p$ with $a_j\le p-1<p\le a_{j+1}$.\vspace{6pt}

However, when considering $\widehat{N}:=(p!^2N_p)_{p\in\NN}=\pi^2(N)$ and so multiplying $N$ (or any log-convex sequence) point-wise with $p!^2$, then one can always ensure \hyperlink{beta1}{$(\beta_1)$} because
$$\forall\;Q\in\NN_{\ge 2}:\;\;\;\liminf_{p\rightarrow\infty}\frac{\widehat{\nu}_{Qp}}{\widehat{\nu}_p}=\liminf_{p\rightarrow\infty}Q^2\frac{\nu_{Qp}}{\nu_p}\ge Q^2>Q.$$
Thus $\widehat{N}$ has \hyperlink{beta1}{$(\beta_1)$} and, as pointed out above, $\widehat{N}$ also violates \eqref{genmg}. Consequently, by Proposition \ref{strangechar} the matrix associated with $\omega_{\widehat{N}}$ also violates both \eqref{rstrange} and \eqref{bstrange}.\vspace{6pt}

In order to guarantee the quasianalyticity part of $(a)$ it is enough to assume that the sequence $(a_j)_j$ satisfies in addition
\begin{equation}\label{sequenceamod}
\forall\;j\ge 1:\;\;\;a_{j+1}\ge\frac{b_j}{j}-a_j.
\end{equation}
This is possible since by \eqref{bproperty} the slope $b_j$ is only depending on $b_1,\dots,b_{j-1}$ and on $a_1,\dots,a_j$. Because $\nu_p=b_j$ for all $p\in\NN$ with $a_j+1\le p\le a_{j+1}$, $j\ge 1$, by \eqref{sequenceamod} we get
$$\sum_{p\ge 1}\frac{1}{\nu_p}=\sum_{j\ge 1}\frac{a_{j+1}-a_j}{b_j}\ge\sum_{j\ge 1}\frac{1}{j}=+\infty.$$

\vspace{6pt}
Finally we turn to $(b)$ and the aim is to modify slightly the construction of $N$ in order to violate not only \eqref{genmg} but also \eqref{equlemma1equ1}. For this we assume that the sequence of slopes $(b_j)_{j\ge 1}$ satisfies
\begin{equation}\label{bpropertynew}
b_j\ge\frac{f(a_j)+j^2(a_j+1)+2j^2(a_j+1)^2}{a_j},\;\;\;\forall\;j\ge 2,
\end{equation}
which is stronger than \eqref{bproperty} above.\vspace{6pt}

Note that $N$ has \eqref{equlemma1equ1} if $\exists\;d\in\NN_{\ge 2}\;\exists\;A,C\ge 1\;\forall\;p\in\NN_{>0}:$
$$\nu_p\le AC^{2p}(N_{dp})^{1/(dp)}\Leftrightarrow dp(f(p)-f(p-1)-\log(A)-2p\log(C))\le f(dp).$$

Let now $d\in\NN_{>0}$, $A, C\ge 1$ be given, arbitrary (large) but from now on fixed. Then for all $l\ge d\ge 2$ the choice $p=a_l+1$ yields
$$(f(p)-f(p-1)-\log(A)-2p\log(C))dp=(b_l-\log(A)-2(a_l+1)\log(C))d(a_l+1),$$
whereas $f(dp)=f(a_l)+b_l((a_l+1)d-a_l)$ (recall the choice in \eqref{aproperty}). Then
\begin{align*}
&(b_l-\log(A)-2(a_l+1)\log(C))d(a_l+1)\le f(a_l)+b_l((a_l+1)d-a_l)
\\&
\Leftrightarrow -\log(A)d(a_l+1)-2(a_l+1)^2d\log(C)\le f(a_l)-a_lb_l
\\&
\Leftrightarrow b_l\le\frac{f(a_l)+\log(A)d(a_l+1)+2\log(C)d(a_l+1)^2}{a_l},
\end{align*}
which leads to a contradiction as $l\rightarrow\infty$ by the modified choice of $b_l$ in \eqref{bpropertynew}. Thus in this case even \eqref{equlemma1equ1} is violated for $N$.\vspace{6pt}

Of course, a combination of $(a)$ and $(b)$ is possible when considering for the constructed sequence in $(b)$ again $\widehat{N}:=(p!^2N_p)_{p\in\NN}$, since again $\widehat{N}$ has to violate \eqref{equlemma1equ1} too, resp. by considering for $(a_j)_j$ again the modified choice \eqref{sequenceamod}.
\qed\enddemo

Theorem \ref{counter1} yields the following consequences when considering the matrix $\mathcal{M}_{\omega_N}$.

\begin{itemize}
\item[$(i)$] Assumptions \eqref{rstrange} and \eqref{bstrange} are in general not guaranteed for weight matrices $\mathcal{M}_{\omega}$ associated with Braun-Meise-Taylor weight functions $\omega$ and hence have to be assumed in addition when required in proofs and arguments, see e.g. in \cite{whitneyextensionweightmatrix} and the additional assumptions \cite[Cor. 9, $(28)$, $(29)$]{whitneyextensionmixedweightfunctionII}, $(28)$ corresponds to \eqref{rstrange} and $(29)$ to \eqref{bstrange}.

\item[$(ii)$] However, in the weight function setting we point out that one may avoid this (technical) problem when using different arguments: For this one should compare \cite[Corollary 9]{whitneyextensionmixedweightfunctionII} with \cite[Thm. 11]{whitneyextensionmixedweightfunctionII} and the more general result \cite[Thm. 4.8]{almostanalytic}.

\item[$(iii)$] Based on this observation a conjecture of the author has been that for weight functions that are concave (resp. equivalent to a concave weight), the associated matrix does always have \eqref{rstrange} and \eqref{bstrange}. In particular this fact should then be true for strong weights since each strong weight is equivalent to a concave one by \cite[Prop. 1.3]{MeiseTaylor88}. However, $(a)$ in Theorem \ref{counter1} shows that this conjecture is not true.

\item[$(iv)$] $(a)$ in this result also shows that both requirements \eqref{rstrange} and \eqref{bstrange} are independent of the notion of (non)quasianalyticity for weight sequences, associated weight functions and their associated matrices, see \cite[Lemma 4.1]{Komatsu73} and \cite[Sect. 4, Cor. 4.8]{testfunctioncharacterization}.
\end{itemize}

We close by summarizing some more facts for $N\in\hyperlink{LCset}{\mathcal{LC}}$ satisfying $(a)$ and/or $(b)$ in Theorem \ref{counter1}.

\begin{lemma}
Let $N\in\hyperlink{LCset}{\mathcal{LC}}$ be a sequence as constructed in Theorem \ref{counter1}. Then we get:
\begin{itemize}
\item[$(i)$] If $N$ is in addition strong non-quasianalytic, then for any $L\in\hyperlink{LCset}{\mathcal{LC}}$ satisfying $L\hyperlink{approx}{\approx}N$ we get $\omega_N\hyperlink{sim}{\sim}\omega_L$ and the matrices $\mathcal{M}_{\omega_N}$ and $\mathcal{M}_{\omega_L}$ are both $R$- and $B$-equivalent.

\item[$(ii)$] If $N$ has in addition $(b)$, then even for each $L\in\hyperlink{LCset}{\mathcal{LC}}$ being equivalent to $N$ the matrix $\mathcal{M}_{\omega_L}$ associated with the weight $\omega_L$ also violates both \eqref{rstrange} and \eqref{bstrange}.
\end{itemize}
\end{lemma}

\demo{Proof}
$(i)$ First, $L\hyperlink{approx}{\approx}N$ and the characterizations obtained in \cite{petzsche} yield \hyperlink{beta1}{$(\beta_1)$} for $L$ as well. Then \cite[Proposition 4.4]{Komatsu73} applied to $N$ and $L$ implies that both weights $\omega_M$ and $\omega_L$ are strong and it is known that each strong weight function is equivalent to a concave weight, more precisely to $\kappa_{\omega}(y):=\int_1^{+\infty}\frac{\omega(yt)}{t^2}dt$, see \cite[Prop. 1.3]{MeiseTaylor88}. In particular both weights are having \hyperlink{om1}{$(\omega_1)$}, too. So \cite[Remark 3.3]{sectorialextensions} yields $\omega_N\hyperlink{sim}{\sim}\omega_L$, see also \cite[Lemma 3.18 $(1)$]{testfunctioncharacterization}.

The fact that $\mathcal{M}_{\omega_N}$ and $\mathcal{M}_{\omega_L}$ are both $R$- and $B$-equivalent follows from $\omega_N\hyperlink{sim}{\sim}\omega_L$ and \cite[Lemma 5.16]{compositionpaper}.\vspace{6pt}

$(ii)$ This follows by Lemma \ref{equlemma} and Proposition \ref{strangechar}.
\qed\enddemo

\appendix

\section{On the generalized strong non-quasianalyticity condition}\label{snqappendix}

By exploiting the formulas for the quotient sequences in the the proof of Proposition \ref{strangechar} we obtain the following statement.

\begin{proposition}\label{strongnonquasilemma}
Let $\omega\in\hyperlink{omset0}{\mathcal{W}_0}$ be given and let $\mathcal{M}_{\omega}=\{W^{(l)}: l>0\}$ be the associated weight matrix. For any $x>0$ and any given $\beta\ge 0$ (fixed) we consider the following growth property:
\begin{equation}\label{condvinthm311}
\exists\;Q\in\NN_{>0}:\;\;\;\liminf_{p\rightarrow+\infty}\frac{\vartheta^{(x)}_{Qp}}{\vartheta^{(x)}_p}>Q^{\beta},
\end{equation}
i.e. \cite[Thm. 3.11 $(v)$]{index} for $W^{(x)}$. In particular, the choice $\beta=0$ yields \eqref{beta3} (condition $(\beta_3)$) and $\beta=1$ yields \hyperlink{beta1}{$(\beta_1)$} for $W^{(x)}$. Then we get:

\begin{itemize}
\item[$(I)$] For any $\beta\ge 0$ the following are equivalent:
\begin{itemize}
\item[$(i)$] There exists $x>0$ such that $W^{(x)}$ satisfies \eqref{condvinthm311}.

\item[$(ii)$] There exists $x>0$ such that for all $c\in\NN_{>0}$ the sequences $W^{(cx)}$ satisfy \eqref{condvinthm311} (uniformly with the same choice for $Q$).
\end{itemize}

\item[$(II)$] Moreover, the following are equivalent:

\begin{itemize}
\item[$(i)$] There exists $x>0$ such that $W^{(x)}$ satisfies \eqref{condvinthm311} with $\beta=0$ (i.e. \eqref{beta3}, $(\beta_3)$).

\item[$(ii)$] There exists $x>0$ such that for all $c\in\NN_{>0}$ the sequences $W^{(cx)}$ and $W^{(x/c)}$ satisfy \eqref{condvinthm311} with $\beta=0$.
\end{itemize}
\end{itemize}
\end{proposition}

\demo{Proof}
The implications $(I)(ii)\Rightarrow(i)$ and $(II)(ii)\Rightarrow(i)$ are trivial (take $c=1$).\vspace{6pt}

$(I)(i)\Rightarrow(ii)$ Let $x>0$ be such that $W^{(x)}$ has \eqref{condvinthm311} with some $Q\in\NN_{>0}$ (in fact one has to take $Q\ge 2$). Then for all $c\in\NN_{>0}$ and $p\in\NN_{>0}$ we obtain (see \eqref{beta1forx1})
$$\vartheta^{(cx)}_p=\frac{W^{(cx)}_p}{W^{(cx)}_{p-1}}=\left(\frac{W^{(x)}_{cp}}{W^{(x)}_{c(p-1)}}\right)^{1/c}=(\vartheta^{(x)}_{c(p-1)+1}\cdots\vartheta^{(x)}_{cp})^{1/c},$$
and for $p=0$ we have clearly the equality $1=1$.

By assumption $\liminf_{p\rightarrow\infty}\frac{\vartheta^{(x)}_{Qp}}{\vartheta^{(y)}_p}>Q^{\beta}$ and by the above identity
$$\frac{\vartheta^{(cx)}_{Qp}}{\vartheta^{(cx)}_p}=\left(\frac{\vartheta^{(x)}_{Qcp-c+1}\cdots\vartheta^{(x)}_{Qcp}}{\vartheta^{(x)}_{cp-c+1}\cdots\vartheta^{(x)}_{cp}}\right)^{1/c}.$$
For all $1\le l\le c$ we have
\begin{equation*}\label{strongnonquasilemmaequ}
Qcp-c+l\ge Q(cp-c+l)\Leftrightarrow Qcp-c+l\ge Qcp-Qc+Ql\Leftrightarrow Q(c-l)\ge(c-l),
\end{equation*}
and so by this and log-convexity for $W^{(x)}$ we obtain
$$\liminf_{p\rightarrow\infty}\frac{\vartheta^{(cx)}_{Qp}}{\vartheta^{(cx)}_p}\ge\prod_{i=1}^c\left(\liminf_{p\rightarrow\infty}\frac{\vartheta^{(x)}_{Qcp-c+i}}{\vartheta^{(x)}_{cp-c+i}}\right)^{1/c}\ge\prod_{i=1}^c\left(\liminf_{p\rightarrow\infty}\frac{\vartheta^{(x)}_{Q(cp-c+i)}}{\vartheta^{(x)}_{cp-c+i}}\right)^{1/c}>Q^{\beta},$$ verifying \eqref{condvinthm311} for $W^{(cx)}$ with the same $Q$.\vspace{6pt}

$(II)(i)\Rightarrow(ii)$ The proof for $W^{(cx)}$ is a particular case of $(I)(i)\Rightarrow(ii)$ and we have only to deal with $W^{(x/c)}$.

By inspecting the proof of \eqref{beta1forxbeur1} we see that the following equality holds true:
$$\vartheta^{(x/c)}_{cp}=\frac{W^{(x/c)}_{cp}}{W^{(x/c)}_{cp-1}}=\frac{W^{(x/c)}_{cp}}{W^{(x/c)}_{cp-c}}\frac{1}{\vartheta^{(x/c)}_{cp-1}}\frac{1}{\vartheta^{(x/c)}_{cp-2}}\cdots\frac{1}{\vartheta^{(x/c)}_{cp-c+1}}=\left(\frac{W^{(x)}_{p}}{W^{(x)}_{p-1}}\right)^c\frac{1}{\vartheta^{(x/c)}_{cp-1}}\frac{1}{\vartheta^{(x/c)}_{cp-2}}\cdots\frac{1}{\vartheta^{(x/c)}_{cp-c+1}}.$$
This implies
$$\forall\;x>0\;\forall\;c\in\NN_{>0}\;\forall\;p\in\NN_{>0}:\;\;\;\vartheta^{(x)}_p=(\vartheta^{(x/c)}_{cp}\vartheta^{(x/c)}_{cp-1}\cdots\vartheta^{(x/c)}_{cp-c+1})^{1/c},$$
and so by involving the log-convexity for $W^{(x/c)}$ and since $Qcp\le 2Qc(p-1)\Leftrightarrow 2\le p$ we can estimate as follows (for any $x>0$, $c,Q\in\NN_{>0}$ and $p\ge 2$):
\begin{equation}\label{strongnonquasilemmaequ1}
\frac{\vartheta^{(x)}_{Qp}}{\vartheta^{(x)}_p}=\left(\frac{\vartheta^{(x/c)}_{Qcp}\vartheta^{(x/c)}_{Qcp-1}\cdots\vartheta^{(x/c)}_{Qcp-c+1}}{\vartheta^{(x/c)}_{cp}\vartheta^{(x/c)}_{cp-1}\cdots\vartheta^{(x/c)}_{cp-c+1}}\right)^{1/c}\le\frac{\vartheta^{(x/c)}_{Qcp}}{\vartheta^{(x/c)}_{cp-c+1}}\le\frac{\vartheta^{(x/c)}_{2Qc(p-1)}}{\vartheta^{(x/c)}_{c(p-1)}}.
\end{equation}
Since by assumption $\liminf_{p\rightarrow+\infty}\frac{\vartheta^{(x)}_{Qp}}{\vartheta^{(x)}_p}>1$ we have verified \eqref{condvinthm311} with $\beta=0$ for $W^{(x/c)}$ when choosing $Q':=2Q$ and restricting in the $\liminf$ to all $q\in\NN$ satisfying $q=c(p-1)$ for some $p\ge 2$.

If now $c(p-1)<q<cp$ for some $p\ge 2$, then by log-convexity for $W^{(x/c)}$ and \eqref{strongnonquasilemmaequ1} applied to $p':=2(p-1)+1=2p-1$ we get
$$\frac{\vartheta^{(x)}_{Q(2p-1)}}{\vartheta^{(x)}_{2p-1}}\le\frac{\vartheta^{(x/c)}_{4Qc(p-1)}}{\vartheta^{(x/c)}_{2c(p-1)}}\le\frac{\vartheta^{(x/c)}_{4Qc(p-1)}}{\vartheta^{(x/c)}_{cp}}\le\frac{\vartheta^{(x/c)}_{4Qq}}{\vartheta^{(x/c)}_{q}},$$
since $2c(p-1)\ge cp\Leftrightarrow p\ge 2$. Altogether, we have shown that
$$\liminf_{p\rightarrow+\infty}\frac{\vartheta^{(x/c)}_{4Qp}}{\vartheta^{(x/c)}_p}>1,$$
i.e. \eqref{condvinthm311} for $W^{(x/c)}$ with $\beta=0$ when taking $Q':=4Q$.
\qed\enddemo

\begin{remark}\label{beta1conjecture}
Concerning this result we point out:

\begin{itemize}
\item[$(i)$] Proposition \ref{strongnonquasilemma} can be applied to $\omega\equiv\omega_M$ provided $M\equiv M^{(1)}\in\hyperlink{LCset}{\mathcal{LC}}$ is satisfying \eqref{condvinthm311}.

\item[$(ii)$] Proposition \ref{strongnonquasilemma} gives a partial answer to the following conjecture of the author: Growth and regularity properties hold for some $W^{(x)}$ if and only if for any $W^{(x)}$, analogously as it is known for \hyperlink{mg}{$(\on{mg})$}, see $(iii)$ in Sect. \ref{classesweightmatrices}. More precisely, we expect that $W^{(x)}$ has \hyperlink{beta1}{$(\beta_1)$} or \eqref{beta3} for some $x>0$ if and only if each $W^{(x)}\in\mathcal{M}_{\omega}$ does so.

However, $(I)$ in Proposition \ref{strongnonquasilemma} is sufficient for the Roumieu-case since the matrix $\{W^{(cx)}: c\in\NN_{>0}\}$ is $R$-equivalent to $\mathcal{M}_{\omega}$ (when $x>0$ being fixed). Analogously, $(II)$ is also sufficient for the Beurling case since $\{W^{(x/c)}: c\in\NN_{>0}\}$ is $B$-equivalent to $\mathcal{M}_{\omega}$ (when $x>0$ being fixed).

\item[$(iii)$] In this context recall that all $W^{(x)}$ have \hyperlink{beta1}{$(\beta_1)$} for the weights $\omega_s(t):=\max\{0,\log(t)^s\}$, $s>1$, see \cite[Prop. 5.14]{whitneyextensionweightmatrix}, and each $\omega_s$ has \eqref{assostrongnq}. Note that \hyperlink{beta1}{$(\beta_1)$} for $M$ implies that $\omega_M$ has \eqref{assostrongnq}, see \cite[Prop. 4.4]{Komatsu73}.

\item[$(iv)$] Finally, concerning $(II)$ in Proposition \ref{strongnonquasilemma} we recall:

By \eqref{goodequivalenceclassic} we get $\omega\hyperlink{sim}{\sim}\omega_{W^{(l)}}$ for each $l>0$. \eqref{condvinthm311} for $\beta=0$ and for $W^{(l)}$ yields property \hyperlink{om1}{$(\omega_1)$} for $\omega_{W^{(l)}}$, see \cite[Lemma 12, $(2)\Rightarrow(4)$]{BonetMeiseMelikhov07}. Hence, if some $W^{(x)}$ has $(II)(i)$, then \hyperlink{om1}{$(\omega_1)$} follows for $\omega$ and for any $\omega_{W^{(l)}}$ since this property is clearly preserved under \hyperlink{sim}{$\sim$}.

Note that the proof of $(II)(i)\Rightarrow(ii)$ shows the following: Assume that $W^{(x)}$ satisfies \eqref{condvinthm311} for some $\beta\ge 0$ with the choice $Q$. Then each $W^{(x/c)}$, $c\in\NN_{>0}$, has \eqref{condvinthm311} with the uniform choice $Q':=4Q$ for $\beta':=\beta(1-\log(4)/\log(4Q))$, hence $0\le\beta'\le\beta$.
\end{itemize}
\end{remark}

By combining this information we close with the following application: Let $M\in\hyperlink{LCset}{\mathcal{LC}}$ be given and satisfying

\begin{itemize}
\item[$(I)$] \hyperlink{beta1}{$(\beta_1)$},

\item[$(II)$] \eqref{genmg},

\item[$(III)$] and finally
\begin{equation}\label{mualmost}
\exists\;A\ge 1\;\forall\;p\in\NN:\;\;\;\mu_{p+1}\le A\mu_p.
\end{equation}
\end{itemize}
Note that $\hyperlink{mg}{(\on{mg})}\Rightarrow\eqref{mualmost}\Rightarrow\hyperlink{dc}{(\on{dc})}$ and in general each implication cannot be reversed, see \cite[Rem. 2.1.36, p. 78]{dissertationjimenez}. The so-called $q$-Gevrey sequences $(q^{p^2})_{p\ge 0}$, $q>1$, satisfy all these requirements.\vspace{6pt}

Let us show that the matrix $\mathcal{M}'_{\omega_M}:=\{M^{(c)}: c\in\NN_{>0}\}$, which is $R$-equivalent to $\mathcal{M}_{\omega_M}$, is ''admissible'' in the notion of \cite[Def. 4.6]{whitneyextensionweightmatrix} (and so $\omega_M$ is an ''admissible'' weight function).

First, by \eqref{beta1forx1} it is immediate to see that each $M^{(c)}$ has \eqref{mualmost} as well since by iterating this property we get $\mu_{cp+i}\le A^c\mu_{cp-c+i}$ for all $p,c\in\NN_{>0}$ and $1\le i\le c$ and so
$$\mu^{(c)}_{p+1}=(\mu_{cp+1}\cdots\mu_{cp+c})^{1/c}\le(A^c\mu_{cp-c+1}\cdots A^c\mu_{cp})^{1/c}=A^c\mu^{(c)}_p.$$

Then \cite[Def. 4.6 $(1)$, $(5)$]{whitneyextensionweightmatrix} for $\mathcal{M}'_{\omega_M}$ hold by definition and Corollary \ref{matrixmgtheoremcor}, whereas \cite[Def. 4.6 $(4)$]{whitneyextensionweightmatrix} follows by $(II)$ and Proposition \ref{strangechar}. The first part of Proposition \ref{strongnonquasilemma} applied to $\beta=1$ and
$(I)$ yield that each $M^{(c)}$, $c\in\NN_{>0}$, has \hyperlink{beta1}{$(\beta_1)$} (i.e. is strongly non-quasianalytic). So $\mathcal{M}'_{\omega_M}$ consists only of strong non-quasianalytic sequences and thus each $M^{(c)}$ is equivalent to its so-called ''descendant'' $S^{(c)}$, see \cite[Sect. 4.1]{whitneyextensionweightmatrix}. Hence \cite[Def. 4.6 $(2)$]{whitneyextensionweightmatrix} is valid, too.

Finally, \cite[Def. 4.6 $(3)$]{whitneyextensionweightmatrix} follows since as shown before each $\mu^{(c)}$ has \eqref{mualmost} and since the equivalence between $M^{(c)}$ and its descendant $S^{(c)}$ even holds on the level of the corresponding quotient sequences, see \cite[Lemma 4.2 $(1)$, $(5)$]{whitneyextensionweightmatrix}.\vspace{6pt}

Summarizing we can apply \cite[Cor. 5.5]{whitneyextensionweightmatrix} as well as the characterizing result \cite[Thm. 5.12]{whitneyextensionweightmatrix} to $\mathcal{M}'_{\omega_M}$.

\bibliographystyle{plain}
\bibliography{Bibliography}

\begin{thebibliography}{10}

\bibitem{Bjorck66}
G.~Bj{\"o}rck.
\newblock Linear partial differential operators and generalized distributions.
\newblock {\em Ark. Mat.}, 6:351--407 (1966), 1966.

\bibitem{BonetBraunMeiseTaylorWhitneyextension}
J.~Bonet, R.~W. Braun, R.~Meise, and B.~A. Taylor.
\newblock Whitney's extension theorem for nonquasianalytic classes of
  ultradifferentiable functions.
\newblock {\em Studia Mathematica}, 99(2):155--184, 1991.

\bibitem{BonetMeiseMelikhov07}
J.~Bonet, R.~Meise, and S.~N. Melikhov.
\newblock A comparison of two different ways to define classes of
  ultradifferentiable functions.
\newblock {\em Bull. Belg. Math. Soc. Simon Stevin}, 14:424--444, 2007.

\bibitem{BraunMeiseTaylor90}
R.~W. Braun, R.~Meise, and B.~A. Taylor.
\newblock Ultradifferentiable functions and {F}ourier analysis.
\newblock {\em Results Math.}, 17(3-4):206--237, 1990.

\bibitem{almostanalytic}
S.~F\"{u}rd\"{o}s, D.~N. Nenning, A.~Rainer, and G.~Schindl.
\newblock Almost analytic extensions of ultradifferentiable functions with
  applications to microlocal {A}nalysis.
\newblock {\em J. Math. Anal. Appl.}, 481(1):123451, 2020.

\bibitem{dissertationjimenez}
J.~Jiménez-Garrido.
\newblock Applications of regular variation and proximate orders to
  ultraholomorphic classes, asymptotic expansions and multisummability.
\newblock PhD Thesis, Universidad de Valladolid, 2018, available online at
  \url{http://uvadoc.uva.es/handle/10324/29501}.

\bibitem{JimenezGarridoSanz}
J.~Jiménez-Garrido and J.~Sanz.
\newblock Strongly regular sequences and proximate orders.
\newblock {\em J. Math. Anal. Appl.}, 438(2):920--945, 2016.

\bibitem{index}
J.~Jiménez-Garrido, J.~Sanz, and G.~Schindl.
\newblock Indices of {O}-regular variation for weight functions and weight
  sequences.
\newblock {\em Rev. R. Acad. Cienc. Exactas Fís. Nat. Ser. A. Mat. RACSAM},
  113(4):3659--3697, 2019.

\bibitem{sectorialextensions}
J.~Jiménez-Garrido, J.~Sanz, and G.~Schindl.
\newblock Sectorial extensions, via {L}aplace transforms, in ultraholomorphic
  classes defined by weight functions.
\newblock {\em Results Math.}, 74(27), 2019.

\bibitem{mixedgrowthindex}
J.~Jiménez-Garrido, J.~Sanz, and G.~Schindl.
\newblock Equality of ultradifferentiable classes by means of indices of mixed
  {O}-regular variation.
\newblock {\em Res. Math.}, 77:art. no. 28, 2022.

\bibitem{Komatsu73}
H.~Komatsu.
\newblock Ultradistributions. {I}. {S}tructure theorems and a characterization.
\newblock {\em J. Fac. Sci. Univ. Tokyo Sect. IA Math.}, 20:25--105, 1973.

\bibitem{KMRc}
A.~Kriegl, P.~W. Michor, and A.~Rainer.
\newblock The convenient setting for non-quasianalytic {D}enjoy--{C}arleman
  differentiable mappings.
\newblock {\em J. Funct. Anal.}, 256:3510--3544, 2009.

\bibitem{mandelbrojtbook}
S.~Mandelbrojt.
\newblock {\em Séries adhérentes, Régularisation des suites, Applications}.
\newblock Gauthier-Villars, Paris, 1952.

\bibitem{matsumoto}
W.~Matsumoto.
\newblock Characterization of the separativity of ultradifferentiable classes.
\newblock {\em J. Math. Kyoto Univ.}, 24(4):667--678, 1984.

\bibitem{matsumotopseudo}
W.~Matsumoto.
\newblock Theory of pseudo-differential operators of ultradifferentiable class.
\newblock {\em J. Math. Kyoto Univ.}, 27(3):453--500, 1987.

\bibitem{MeiseTaylor88}
R.~Meise and B.~A. Taylor.
\newblock Whitney's extension theorem for ultradifferentiable functions of
  {B}eurling type.
\newblock {\em Ark. Mat.}, 26(2):265--287, 1988.

\bibitem{petzsche}
H.-J. Petzsche.
\newblock On {E}. {B}orel's theorem.
\newblock {\em Math. Ann.}, 282(2):299--313, 1988.

\bibitem{PetzscheVogt}
H.-J. Petzsche and D.~Vogt.
\newblock Almost analytic extension of ultradifferentiable functions and the
  boundary values of holomorphic functions.
\newblock {\em Math. Ann.}, 267:17--35, 1984.

\bibitem{compositionpaper}
A.~Rainer and G.~Schindl.
\newblock Composition in ultradifferentiable classes.
\newblock {\em Studia Math.}, 224(2):97--131, 2014.

\bibitem{whitneyextensionweightmatrix}
A.~Rainer and G.~Schindl.
\newblock Extension of {W}hitney jets of controlled growth.
\newblock {\em Math. Nachr.}, 290(14--15):2356--2374, 2017.

\bibitem{whitneyextensionmixedweightfunctionII}
A.~Rainer and G.~Schindl.
\newblock On the extension of {W}hitney ultrajets, {I}{I}.
\newblock {\em Studia Math.}, 250(3):283--295, 2020.

\bibitem{dissertation}
G.~Schindl.
\newblock Exponential laws for classes of {D}enjoy-{C}arleman-differentiable
  mappings, 2014.
\newblock PhD Thesis, Universität Wien, available online at
  \url{http://othes.univie.ac.at/32755/1/2014-01-26_0304518.pdf}.

\bibitem{testfunctioncharacterization}
G.~Schindl.
\newblock Characterization of ultradifferentiable test functions defined by
  weight matrices in terms of their {F}ourier transform.
\newblock {\em Note di Matematica}, 36(2):1--35, 2016.

\bibitem{subaddlike}
G.~Schindl.
\newblock On subadditivity-like conditions for associated weight functions.
\newblock {\em Bull. Belg. Math. Soc. Simon Stevin}, 28(3):399--427, 2022.

\end{thebibliography}

\end{document}